\pdfoutput=1
\RequirePackage{ifpdf}
\ifpdf 
\documentclass[pdftex]{sigma}
\else
\documentclass{sigma}
\fi

\newtheorem{Theorem}{Theorem}[section]

\newtheorem{thm}[Theorem]{Theorem}

\newtheorem{lem}[Theorem]{Lemma}

\theoremstyle{definition}
\newtheorem{rem}[Theorem]{Remark}

\newtheorem{Example}[Theorem]{Example}

\numberwithin{equation}{section}

\newcommand{\R}{\mathbb{R}}
\def\C{\mathbb{C}}
\newcommand{\Q}{\mathbb{Q}}
\newcommand{\Z}{\mathbb{Z}}
\newcommand{\g}{\mathfrak{g}}
\newcommand{\ga}{\gamma}
\newcommand{\nc}{\newcommand}
\nc{\on}{\operatorname}
\nc{\la}{\lambda}
\nc{\wh}{\widehat}
\nc{\wt}{\widetilde}
\nc{\ghat}{\wh{\g}}
\nc{\mc}{\mathcal}
\nc{\pa}{\partial}
\nc{\ppart}{(\!(t)\!)}
\nc{\zpart}{(\!(z)\!)}
\nc{\n}{{\mathfrak n}}
\nc{\ol}{\overline}
\nc{\mb}{\mathbf}
\nc{\h}{{\mathfrak h}}
\nc{\pone}{{\mathbb P}^1}
\nc{\bs}{\backslash}
\nc{\al}{\alpha}
\nc{\gt}{{\mathfrak g}'}
\nc{\Bun}{\on{Bun}}
\def\neg{}
\nc{\ka}{\kappa}

\def\LG{{}^L\neg G}
\def\LT{{}^L\neg T}
\nc{\OO}{\mathcal O}
\nc{\Fq}{{\mathbb F}_q}
\nc{\cmu}{\check\mu}
\nc{\hl}{h_{\ell,x}}
\nc{\hr}{h_{r,x}}
\nc{\M}{{\mc M}}
\nc{\K}{{\mc K}}
\nc{\Pic}{\on{Pic}}

\begin{document}
\allowdisplaybreaks

\newcommand{\arXivNumber}{1812.08160}

\renewcommand{\thefootnote}{}

\renewcommand{\PaperNumber}{042}

\FirstPageHeading

\ShortArticleName{Is There an Analytic Theory of Automorphic Functions for Complex Algebraic Curves?}

\ArticleName{Is There an Analytic Theory of Automorphic \\ Functions for Complex Algebraic Curves?\footnote{This paper is a~contribution to the Special Issue on Algebra, Topology, and Dynamics in Interaction in honor of Dmitry Fuchs. The full collection is available at \href{https://www.emis.de/journals/SIGMA/Fuchs.html}{https://www.emis.de/journals/SIGMA/Fuchs.html}}}

\Author{Edward FRENKEL}
\AuthorNameForHeading{E.~Frenkel}
\Address{Department of Mathematics, University of California,
 Berkeley, CA 94720, USA}
\Email{\href{mailto:frenkel@math.berkeley.edu}{frenkel@math.berkeley.edu}}

\ArticleDates{Received September 30, 2019, in final form April 27, 2020; Published online May 16, 2020}

\Abstract{The geometric Langlands correspondence for complex algebraic
  curves differs from the original Langlands correspondence for number
  fields in that it is formulated in terms of sheaves rather than
  functions (in the intermediate case of curves over finite fields,
  both formulations are possible). In a recent preprint, Robert
  Langlands made a proposal for developing an analytic theory of
  automorphic forms on the moduli space of $G$-bundles on a complex
  algebraic curve. Langlands envisioned these forms as eigenfunctions
  of some analogues of Hecke operators. In these notes I show that if
  $G$ is an abelian group then there are well-defined Hecke operators,
  and I give a complete description of their eigenfunctions and
  eigenvalues. For non-abelian $G$, Hecke operators involve
  integration, which presents some difficulties. However, there is an
  alternative approach to developing an analytic theory of automorphic
  forms, based on the existence of a large commutative algebra of
  global differential operators acting on half-densities on the moduli
  stack of $G$-bundles. This approach (which implements some ideas of
  Joerg Teschner) is outlined here, as a preview of a joint work with
  Pavel Etingof and David Kazhdan.}

\Keywords{Langlands Program; automorphic function; complex algebraic
  curve; principal $G$-bundle; Jacobian variety; differential
  operator; oper}

\Classification{14D24; 17B67; 22E57}

\renewcommand{\thefootnote}{\arabic{footnote}}
\setcounter{footnote}{0}

\rightline{\em To my teacher Dmitry Borisovich Fuchs on his 80th birthday}

\section{Introduction}

{\bf 1.1.}~The foundations of the Langlands Program were laid by
Robert Langlands in the late 1960s~\cite{L}. Originally, these
ideas were applied in two realms: that of number fields, i.e., finite
extensions of the field $\Q$ of rational numbers, and that of function
fields, where by a~function field one understands the field of
rational functions on a smooth projective curve over a finite field~$\Fq$. In both cases, the objects of interest are {\em automorphic
 forms}, which are, roughly speaking, functions on the quotient of
the form $G(F)\bs G({\mathbb A}_F)/K$. Here $F$ is a number field or a
function field, $G$ is a reductive
algebraic group over~$F$, ${\mathbb A}_F$ is the ring of adeles of~$F$, and~$K$
is a compact subgroup of~$G({\mathbb A}_F)$. There is a family of
mutually commuting {\em Hecke operators} acting on this space of
functions, and one wishes to describe the common eigenfunctions of
these operators as well as their eigenvalues. The idea is that those
eigenvalues can be packaged as the ``Langlands parameters'' which can
be described in terms of homomorphisms from a group closely related to
the Galois group of~$F$ to the Langlands dual group~$\LG$ associated
to~$G$, and perhaps some additional data.

To be more specific, let $F$ be the field of rational functions on a
curve $X$ over $\Fq$ and $G={\rm GL}_n$. Let us further restrict ourselves
to the unramified case, so that $K$ is the maximal compact subgroup
$K={\rm GL}_n(\OO_F)$, where $\OO_F \subset {\mathbb A}_F$ is the ring of
integer adeles. In this case, a theorem of V.~Drinfeld \cite{Dr2,Dr1,Dr2-2,Dr2-3}
for $n=2$ and L.~Lafforgue \cite{Laf} for $n>2$ states that (if we
impose the so-called cuspidality condition and place a restriction on
the action of the center of ${\rm GL}_n$) the Hecke eigenfunctions on
${\rm GL}_n(F)\bs {\rm GL}_n({\mathbb A}_F)/{\rm GL}_n(\OO_F)$ are in one-to-one
correspondence with $n$-dimensional irreducible unramified
representations of the Galois group of $F$ (with a matching
restriction on its determinant).

{\bf 1.2.}~Number fields and function fields for curves over $\Fq$
are two ``languages'' in Andr\'e Weil's famous trilingual ``Rosetta
stone'' \cite{Weil}, the third language being the theory of algebraic
curves over the field $\C$ of complex numbers. Hence it is tempting to
build an analogue of the Langlands correspondence in the setting of a
complex curve $X$. Such a theory has indeed been developed starting
from the mid-1980s, initially by V.~Drinfeld \cite{Dr1} and G.~Laumon
\cite{Laumon} (and relying on the ideas of an earlier work of P.~Deligne),
then by A.~Beilinson and V.~Drinfeld \cite{BD}, and subsequently by
many others. See, for example, the surveys \cite{F:rev,Ga:bourbaki}
for more details. However, this theory, dubbed ``geometric Langlands
Program'', is quite different from the Langlands Program in its
original formulation for number fields and function fields.

The most striking difference is that in the geometric theory the
vector space of automorphic functions on the double quotient $G(F)\bs
G({\mathbb A}_F)/K$ is replaced by a (derived) category of sheaves on
an algebraic stack whose set of $\C$-points is this quotient. For
example, in the unramified case $K=G(\OO_F)$, this is the moduli stack~$\Bun_G$ of principal $G$-bundles on our complex curve~$X$. Instead of the Hecke operators of the classical theory, which act on functions,
we then have {\em Hecke functors} acting on suitable categories of
sheaves, and instead of Hecke eigenfunctions we have {\em Hecke eigensheaves}.

For example, in the unramified case a Hecke eigensheaf ${\mc F}$ is a
sheaf on $\Bun_G$ (more precisely, an object in the category of
$D$-modules on $\Bun_G$, or the category of perverse sheaves on
$\Bun_G$) with the property that its images under the Hecke functors
are isomorphic to ${\mc F}$ itself, tensored with a vector space (this
is the categorical analogue of the statement that under the action of
the Hecke operators eigenfunctions are multiplied by
scalars). Furthermore, since the Hecke functors (just like the Hecke
operators acting on functions) are parametrized by closed points of~$X$, a~Hecke eigensheaf actually yields a family of vector spaces parametrized by points of~$X$. We then impose an additional
requirement that these vector spaces be stalks of a local system on~$X$ for the Langlands dual group $\LG$ (taken in the representation of~$\LG$ corresponding to the Hecke functor under consideration). This
neat formulation enables us to directly link Hecke eigensheaves and
(equivalence classes of) $\LG$-local systems on~$X$, which are the
same as (equivalence classes of) homomorphisms from the fundamental group $\pi_1(X,p_0)$ of $X$ to $\LG$.

\looseness=-1 This makes sense from the point of view of Weil's Rosetta stone,
because the fundamental group can be seen as a geometric analogue of
the unramified quotient of the Galois group of a function field. We
note that for $G={\rm GL}_n$, in the unramified case, the Hecke eigensheaves
have been constructed in~\cite{Dr1} for $n=2$ and in
\cite{FGV,Ga1} for $n>2$. More precisely, the following theorem has
been proved: for any irreducible rank~$n$ local system ${\mc E}$ on
$X$, there exists a Hecke eigensheaf on $\Bun_{{\rm GL}_n}$ whose
``eigenvalues'' correspond to~${\mc E}$.\footnote{Furthermore, these
 Hecke eigensheaves are irreducible on each connected component of
 $\Bun_{{\rm GL}_n}$.} Many results of that nature have been obtained for
other groups as well. For example, in \cite{BD} Hecke eigensheaves on
$\Bun_G$ were constructed for all $\LG$-local systems having the
structure of an $\LG$-{\em oper} (these local systems form a
Lagrangian subspace in the moduli of all $\LG$-local
systems). Furthermore, a more satisfying categorical version of the
geometric Langlands correspondence has been proposed by A.~Beilinson
and V.~Drinfeld and developed further in the works of D.~Arinkin and
D.~Gaitsgory \cite{AG,Ga:ast} (see~\cite{Ga:bourbaki} for a survey).

\looseness=1 To summarize, the salient difference between the original formulation
of the Langlands Program (for number fields and function fields of
curves over $\Fq$) and the geometric formulation is that the former is
concerned with functions and the latter is concerned with
sheaves. What makes this geometric formulation appealing is that in
the intermediate case -- that of curves over~$\Fq$~-- which serves as
a kind of a bridge in the Rosetta stone between the number field case
and the case of curves over~$\C$, both function-theoretic and
sheaf-theoretic formulations make sense. Moreover, it is quite common
that the same geometric construction works for curves over $\Fq$ and
$\C$. For example, essentially the same construction produces Hecke
eigensheaves on $\Bun_{{\rm GL}_n}$ for an irreducible rank $n$ local system
on a curve over $\Fq$ and over $\C$ \cite{Dr1,FGV,Ga1}.\footnote{The
 term ``local system'' has different meanings in the two cases: it is
 an $\ell$-adic sheaf in the first case and a bundle with a flat
 connection in the second case, but what we do with these local
 systems to construct Hecke eigensheaves (in the appropriate
 categories of sheaves) is essentially the same in both cases.}

Furthermore, in the realm of curves over $\Fq$, the function-theoretic
and sheaf-theoretic formulations are connected to each other by
Alexander Grothendieck's ``functions-sheaves dictionary''. This
dictionary assigns to a ($\ell$-adic) sheaf~${\mc F}$ on a~variety (or
an algebraic stack)~$V$ over~$\Fq$, a~function on the set of closed
points of~$V$ whose value at a given closed point~$v$ is the
alternating sum of the traces of the Frobenius (a~generator of the
Galois group of the residue field of~$v$) on the stalk cohomologies of
${\mc F}$ at $v$ (see \cite[Section~1.2]{Laumon:const} or \cite[Section~3.3]{F:rev},
for details). Thus, for curves over $\Fq$ the geometric formulation of the Langlands Program may be viewed as a~{\em refinement} of the original formulation: the goal is to produce, for
each $\LG$-local system on~$X$, the corresponding Hecke eigensheaf on~$\Bun_G$, but at the end of the day we can always go back to the more
familiar Hecke eigenfunctions by taking the traces of the Frobenius on
the stalks of the Hecke eigensheaf at the $\Fq$-points of
$\Bun_G$. Thus, the function-theoretic and the sheaf-theoretic
formulations go hand-in-hand for curves over $\Fq$.

{\bf 1.3.}~In the case of curves over $\C$ there is
no Frobenius, and hence no direct way to get functions out of Hecke
eigensheaves on $\Bun_G$. However, since a Hecke eigensheaf is a
$D$-module on $\Bun_G$, we could view its sections as analogues of
automorphic functions of the analytic theory. The problem is that for
non-abelian $G$, these $D$-modules -- and hence their sections -- are
known to have complicated singularities and monodromies. Outside of
the singularity locus, a Hecke eigensheaf is a holomorphic vector
bundle with a holomorphic flat connection, but its horizontal sections
have non-trivial monodromies along the closed paths going around
various components of the singularity locus (and in general there are
non-trivial monodromies along other closed paths as well). So instead
of functions we get multi-valued sections of a vector bundle. On top
of that, in the non-abelian case the rank of this vector bundle grows
exponentially as a~function of the genus of~$X$, and furthermore, the
components of the singularity locus have a rather complicated
structure. Therefore in the non-abelian case, as the genus of~$X$
grows, it becomes increasingly difficult to study these horizontal
sections. For this reason, it is the $D$-modules themselves, rather
than their sections, that are traditionally viewed as more meaningful
objects of study, and that's why in the geometric formulation of the
Langlands Program for curves over~$\C$, we focus on these $D$-modules
rather than their multi-valued sections. Thus, the geometric theory in
the case of complex curves becomes inherently sheaf-theoretic.

{\bf 1.4.}~In a recent preprint \cite{L:analyt}, Robert Langlands made a proposal
for developing an analytic theory of automorphic functions for complex
algebraic curves. He mostly considered the case that~$X$ is an
elliptic curve and~$G$ is ${\rm GL}_1$ or~${\rm GL}_2$. His proposal can be
summarized as follows: (1)~He assumed that one can define a
commutative algebra of Hecke operators acting on a~particular space of
$L^2$ functions on $\Bun_G$ (he only gave a definition of these when
$X$ is an elliptic curve and $G={\rm GL}_2$). (2)~He assumed
that the Satake isomorphism of the theory
over $\Fq$ would also hold over $\C$ and that each point $\sigma$ in
the joint spectrum of these Hecke operators would give rise to a
function $f_\sigma$ on the curve $X$ with values in the space of
semi-simple conjugacy classes of the maximal compact subgroup $\LG_c$
of $\LG$. (3)~He proposed that each function $f_\sigma$ could be
expressed in terms of the holonomies of a Yang--Mills connection
$\nabla_\sigma$ on an $\LG_c$-bundle on $X$. (4)~Atiyah and Bott have
shown in~\cite{AB} that to a Yang--Mills connection $\nabla$ one can
associate a homomorphism $\rho(\nabla)$ from a central extension~$\wh\pi_1(X)$ of the fundamental group~$\pi_1(X)$ of~$X$ to~$\LG_c$. Langlands proposed that the resulting map $\sigma \mapsto
\rho(\nabla_\sigma)$ would give rise to a bijection between the
spectrum of the Hecke operators and the set of equivalence classes of
homomorphisms $\wh\pi_1(X) \to \LG_c$ satisfying a certain finiteness
condition.

{\bf 1.5.}~In this paper, I discuss this proposal. Consider first the case of
${\rm GL}_1$.

In this case, the Picard variety of a complex curve $X$ plays the role
of $\Bun_{{\rm GL}_1}$ (see Section~\ref{anell}). It carries a natural
integration measure using which one can define the Hilbert space of
$L^2$ functions. The Hecke operators are rather simple in the case of
${\rm GL}_1$ (as well as an arbitrary torus): they are pull-backs of
functions under natural maps. Therefore no integration is needed to
define an action of the commutative algebra of Hecke operators on this
Hilbert space. The question of finding their eigenfunctions and
eigenvalues is well-posed.

I give a complete answer to this question in Section \ref{abelian}:
first for elliptic curves in Sections \ref{anell} and \ref{genell} and
then for curves of an arbitrary genus in Section \ref{higher}. In
Section~\ref{torus}, I generalize these results to the case of an
arbitrary torus $T$ instead of ${\rm GL}_1$. In particular, I show that
Hecke eigenfunctions are labeled by $H^1(X,\Lambda^*(T))$, the first
cohomology group of $X$ with coefficients in the lattice of
cocharacters of $T$, and give an explicit formula for the
corresponding eigenvalues. The construction uses the Abel--Jacobi map.

The results presented in Section \ref{abelian} agree with parts (1),
(2), and (3) of Langlands' proposal in the case of ${\rm GL}_1$. However,
the results of Section \ref{abelian} are not in agreement with part
(4) of the proposal. Indeed, each point $\sigma$
in the spectrum of the Hecke operators in the case of ${\rm GL}_1$ gives
rise to a function $f_\sigma$ on~$X$ with values in~$U(1) \subset
\C^\times$ and it is possible to write this function~$f_\sigma$ as the
holonomy of a flat unitary connection $\nabla_\sigma$ on a line bundle
on~$X$. This is shown in Section~\ref{fund} for elliptic curves and in
Section~\ref{higher} for general curves. However, and this is a key
point, each of these connections necessarily gives rises to the {\em
 trivial monodromy representation} of the fundamental group
$\pi_1(X)$. Indeed, by construction, $f_\sigma$ is a {\em
 single-valued} function on~$X$, and it is a horizontal section of
the connection $\nabla_\sigma$. Therefore the
connection~$\nabla_\sigma$ has trivial monodromy. Thus, the map in
part~(4) sends each $\sigma$ to
the trivial representation of $\pi_1(X)$.

{\bf 1.6.}~Now consider the case of ${\rm GL}_2$. Unlike the abelian
case, in order to define Hecke operators for non-abelian groups, one
cannot avoid integration. Therefore one needs to define the pertinent
integration measures. In the classical setting, over $\Fq$, the group
$G(\Fq\ppart)$ is locally compact and therefore carries a Haar
measure. Using this Haar measure, one then defines the measures of
integration pertinent to the Hecke operators. In contrast, the group
$G(\C\ppart)$ is not locally compact, and therefore it does {\em not}
carry a Haar measure, which is only defined for locally compact
groups. Therefore, the standard definition of the measure for curves
over $\Fq$ does not directly generalize to the case of curves over
$\C$, as explained in Section~\ref{nonabelian}.

In~\cite{L:analyt} an attempt is made to explicitly define Hecke
operators acting on a particular version of an $L^2$ space of
$\Bun_{{\rm GL}_2}$ of an elliptic curve. Alas, there are serious issues
with this proposal (see Section~\ref{attempt ell}).

{\bf 1.7.}~There is, however, another possibility: rather than
looking for the eigenfunctions of Hecke operators, one can look for
the eigenfunctions of global differential operators on $\Bun_G$. These
eigenfunctions and the corresponding eigenvalues have been recently
studied for $G={\rm SL}_2$ by Joerg Teschner~\cite{Teschner}. In a joint work with Pavel Etingof and
David Kazhdan~\cite{EFK}, we propose a canonical self-adjoint
extension of the algebra of these differential operators and study the
corresponding spectral problem. I discuss this proposal in Section~\ref{alt}.

According to a theorem of Beilinson and Drinfeld~\cite{BD}, there is a
large commutative algebra of global holomorphic differential operators
acting on sections of a square root $K^{1/2}$ of the canonical line
bundle $K$ on~$\Bun_G$ (this square root always exists, and is unique
if $G$ is simply-connected~\cite{BD}). The complex conjugates of these
differential operators are anti-holomorphic and act on sections of the
complex conjugate line bundle $\ol{K}^{1/2}$ on~$\Bun_G$. The tensor
product of these two algebras is a commutative algebra acting on
sections of the line bundle $K^{1/2} \otimes \ol{K}^{1/2}$ which we
refer to as the bundle of half-densities on $\Bun_G$.

The space of compactly supported sections of the line bundle $K^{1/2}
\otimes \ol{K}^{1/2}$ on $\Bun_G$ (or rather, on its open dense
subspace of stable $G$-bundles, provided that one exists) has a
natural Hermitian inner product. Taking the completion of this space,
we obtain a Hilbert space. Our differential operators are unbounded
linear operators on this Hilbert space. We can ask whether these
operators have natural self-adjoint extensions and if so, what are
their joint eigenfunctions and eigenvalues. In Section \ref{global},
as a preview of \cite{EFK}, I give some more details on this
construction. I then explain what happens in the abelian case of
$G={\rm GL}_1$ in Section \ref{diff GL1}.

\looseness=1 In the case of ${\rm GL}_1$, the global differential operators are
polynomials in the shift vector fields, holomorphic and
anti-holomorphic, on the neutral component $\Pic^0(X)$ of the Picard
variety of a complex curve $X$. These operators commute with each other (and with the Hecke
operators, which are available in the abelian case), and their joint
eigenfunctions are the standard Fourier harmonics on $\Pic^0(X)$. What
about the eigenvalues? The spectrum of the commutative algebra of
global holomorphic differential operators on $\Pic^0(X)$ can be
identified with the space of holomorphic connections on the trivial
line bundle on $X$. Hence every eigenvalue of this algebra can be
encoded by a~point in this space. It turns out that the points
corresponding to the eigenvalues of this algebra on the space of $L^2$
functions on $\Pic^0(X)$ are precisely those holomorphic connections on
the trivial line bundle on $X$ that give rise to the homomorphisms
$\pi_1(X,p_0) \to \C^\times$ with image in $\R^\times \subset
\C^\times$. In other words, these are the connections with monodromy
in the {\em split real form}~${\rm GL}_1(\R)$ of~${\rm GL}_1(\C)$. This dovetails
nicely with the conjecture of Teschner~\cite{Teschner} in the case of
$G={\rm SL}_2$. We expect an analogous statement to hold for a general
reductive group~$G$, see~\cite{EFK}.

\looseness=1 Suppose for simplicity that $G$ is simply-connected. Then, according
to a theorem of Beilinson and Drinfeld~\cite{BD}, the spectrum of the
algebra of global holomorphic differential operators on $\Bun_G$ is
canonically identified with the space of $\LG$-opers on~$X$. If
$G={\rm SL}_2$, then $\LG={\rm PGL}_2$ and ${\rm PGL}_2$-opers are the same as
projective connections. Teschner \cite{Teschner} proposed that in this
case, the eigenvalues correspond to the projective connections with
monodromy taking values in the split real form ${\rm PGL}_2(\R)$ of
${\rm PGL}_2(\C)$ (up to conjugation by an element of ${\rm PGL}_2(\C)$). Such
projective connections have been described by
W.M.~Goldman~\cite{Goldman}. For general $G$, we expect that the
joint eigenvalues of
the global holomorphic differential operators on $\Bun_G$
correspond to those $\LG$-opers that have monodromy taking values in
the split real form of $\LG$ (up to conjugation). If so, then the
spectra of the global differential operators on $\Bun_G$ can be
described by analogues of the Langlands parameters of the classical
theory: namely, certain homomorphisms from the fundamental group of~$X$ to the Langlands dual group~$\LG$. A~somewhat surprising element is that the homomorphisms that appear here are the ones whose image is
in the split real form of~$\LG$ (rather than the compact form). More
details will appear in~\cite{EFK}.

{\bf 1.8.}~Thus, there is a rich analytic theory of joint
eigenfunctions and eigenvalues of the global differential operators
acting on half-densities on $\Bun_G$. This theory can be viewed as
an analytic theory of automorphic functions for complex curves. So,
Langlands was right to insist that an analytic theory exists, and he
deserves a lot of credit for trying to construct it.

This raises the next question: what is the connection between this
analytic theory and the geometric theory?

Valuable insights into this question may be gleaned from
two-dimensional conformal field theory (CFT). In CFT, one has two
types of correlation functions. The first type is chiral correlation
functions, also known as conformal blocks. They form a vector space
for fixed values of the parameters of the CFT. Hence we obtain a
vector bundle of conformal blocks on the space of parameters. In
addition, the data of conformal field theory give rise to a
projectively flat connection on this bundle. The conformal blocks are
{\em multi-valued} horizontal sections of this bundle. The second type
is the ``true'' correlation functions. They can be expressed as
sesquilinear combinations of conformal blocks and their complex
conjugates (anti-conformal blocks), chosen so that the combination is
a {\em single-valued} function of the parameters (see, e.g.,
\cite[Lecture~4]{Gaw}).\footnote{As a useful analogy, consider the
 exponentials of harmonic functions, which may be written as products
 of holomorphic and anti-holomorphic functions.}

Now, the Hecke eigensheaves on $\Bun_G$ constructed in \cite{BD} may
be viewed as sheaves of conformal blocks of a certain two-dimensional
conformal field theory, see \cite{F:rev}. Away from a singularity
locus, these sheaves are vector bundles with a flat connection, and
conformal blocks are their multi-valued horizontal sections (see
Section~1.3 above). It turns out that in some cases there
exist linear combinations of products of these conformal blocks and
their complex conjugates which give rise to single-valued functions on
$\Bun_G$. These functions are precisely the automorphic forms of the
analytic theory. In other words, the objects of the analytic theory of
automorphic forms on $\Bun_G$ can be constructed from the objects of
the geometric Langlands theory in roughly the same way as the
correlation functions of CFT are constructed from conformal
blocks. This was predicted in \cite{F:MSRI} and \cite{Teschner}. An
important difference with the CFT is that whereas in CFT the monodromy
of conformal blocks is typically unitary, here we expect the monodromy
to be in a split real group.

\section{The abelian case} \label{abelian}

\subsection{The case of an elliptic curve} \label{anell}

Let's start with the case of an elliptic curve $E_\tau$ with complex
parameter $\tau$. Let's choose, once and for all, a reference point
$p_0$ on this curve. Then we can identify it with
\begin{gather} \label{Etau}
E_\tau \simeq \C/(\Z + \Z\tau).
\end{gather}

Next, consider the Picard variety $\Pic(E_\tau)$ of $E_\tau$. This is
the (fine) moduli space of line bundles on $E_\tau$ (note that the
corresponding moduli stack $\Bun_{{\rm GL}_1}(E_\tau)$ of line bundles on
$E_\tau$ is the quotient of $\Pic(E_\tau)$ by the trivial action of
the multiplicative group ${\mathbb G}_m={\rm GL}_1$, which is the group of
automorphisms of every line bundle on $E_\tau$). It is a disjoint
union of connected components $\Pic^d(E_\tau)$ corresponding to line
bundles of degree $d$. Using the reference point $p_0$, we can
identity $\Pic^d(E_\tau)$ with $\Pic^0(E_\tau)$ by sending a line
bundle ${\mc L}$ of degree $d$ to ${\mc L}(-d \cdot
p_0)$. Furthermore, we can identify the degree~$0$ component
$\Pic^0(E_\tau)$, which is the Jacobian variety of~$E_\tau$, with
$E_\tau$ itself using the Abel--Jacobi map; namely, we map a point $p
\in E_\tau$ to the degree~0 line bundle~$\OO(p-p_0)$.

Now we define the Hecke operators $H_p$. They are labeled by points
$p$ of the curve $E_\tau$. The operator $H_p$ is the pull-back of
functions with respect to the geometric map
\begin{align} \label{origT}
T_p\colon \ \Pic^d(E_\tau) &\to \Pic^{d+1}(E_\tau), \\
\notag {\mc L} & \mapsto {\mc L}(p).
\end{align}
These operators commute with each other.

Formula \eqref{origT} implies that if $f$ is a joint eigenfunction of
the Hecke operators $H_p$, $p \in E_\tau$, on~$\Pic(E_\tau)$, then its
restriction $f_0$ to the connected component~$\Pic^0(E_\tau)$ is an
eigenfunction of the operators
\begin{gather*}
_{p_0}H_p = H_{p_0}^{-1} H_p,
\end{gather*}
where $p_0$ is our reference point.

Conversely, given an eigenfunction $f_0$ of $_{p_0}H_p$, $p \in X$, on
$\Pic^0(X)$ and $\mu_{p_0} \in \C^\times$, there is a~unique extension
of~$f_0$ to an eigenfunction~$f$ of~$H_p$, $p \in X$, such that the
eigenvalue of $H_{p_0}$ on $f$ is equal to~$\mu_{p_0}$. Namely, any
line bundle ${\mc L}$ of degree~$d$ may be represented uniquely as~${\mc L}_0(d \cdot p_0)$, where ${\mc L}_0$ is a line bundle of degree~$0$. We then set
\begin{gather*}
f({\mc L}) = (\mu_{p_0})^d \cdot f_0({\mc L}_0).
\end{gather*}
By construction, the eigenvalue $\mu_p$ of $H_p$ on $f$ is then equal
to $\la_p \cdot \mu_{p_0}$, where~$\la_p$ is the eigenvalue of~$_{p_0}
H_p$ on~$f_0$ (note that since $_{p_0} H_{p_0} = \on{Id}$, the
eigenvalue $\la_{p_0}$ is always equal to $1$).

Therefore, from now on we will consider the eigenproblem for the
operators $_{p_0}H_p$ acting on the space $L^2\big(\Pic^0(E_\tau)\big)$ of
$L^2$-functions on $\Pic^0(E_\tau)$. Here, we define
$L^2\big(\Pic^0(E_\tau)\big)$ as $L^2(E_\tau)$ (with respect to the measure on
$E_\tau$ induced by the translation-invariant measure on $\C$ via the
isomorphism~\eqref{Etau}) using the above isomorphism between
$\Pic^0(E_\tau)$ and~$E_\tau$. The Hecke operator~$_{p_0}H_p$ acting
on $L^2(E_\tau)$ is given by the formula
\begin{gather} \label{Hpdef}
({}_{p_0}H_p \cdot f)(q) = f(q+p).
\end{gather}
In other words, it is simply the pull-back under the shift by $p$ with
respect to the (additive) abelian group structure on $E_\tau$, which
can be described explicitly using the isomorphism~\eqref{Etau}. The
subscript $p_0$ in $_{p_0}H_p$ serves as a reminder that this operator
depends on the choice of the reference point~$p_0$.

Now we would like to describe the joint eigenfunctions and eigenvalues
of the operators $_{p_0} H_p$ on $L^2(E_\tau)$.

To be even more concrete, let's start with the case $\tau={\rm i}$, so
$E_\tau=E_{\rm i}$ which is identified with $\C/(\Z+\Z {\rm i})$ as above. Thus, we
have a measure-preserving isomorphism between $E_{\rm i}$ and the product of
two circles $(\R/\Z) \times (\R/\Z)$ corresponding to the real and
imaginary parts of $z=x+{\rm i}y$. The space of~$L^2$ functions on the curve
$E_{\rm i}$ is therefore the completed tensor product of two copies of
$L^2(\R/\Z)$, and so it has the standard orthogonal Fourier basis:
\begin{gather} \label{fmn}
f_{m,n}(x,y) = {\rm e}^{2\pi{\rm i} mx} \cdot {\rm e}^{2\pi{\rm i} ny}, \qquad m,n \in \Z.
\end{gather}
Let us write $p = x_p+y_p{\rm i} \in E_{\rm i}$, with $x_p,y_p \in [0,1)$. The
operator $_{p_0}H_p$ corresponds to the shift of $z$ by $p$ (with
respect to the abelian group structure on $E_{\rm i}$):
\begin{gather*} 
({}_{p_0}H_p \cdot f)(x,y) = f(x+x_p,y+y_p), \qquad f \in L^2(E_{\rm i}).
\end{gather*}

It might be instructive to consider first the one-dimensional analogue
of this picture, in which we have $L^2\big(S^1\big)$, where $S_1=\C/\Z$ with
coordinate $\phi$. Then the role of the family $\{ _{p_0} H_p \}_{p\in
 E_{\rm i}}$ is played by the family $\{ H'_{\al} \}_{\al \in S^1}$ acting
by shifts:
\begin{gather*} 
(H'_{\al} \cdot f)(x) = f(\phi+\al), \qquad f \in L^2\big(S^1\big).
\end{gather*}
Then the Fourier harmonics $f_n(x) = {\rm e}^{2\pi{\rm i} n\phi}$ form an
orthogonal eigenbasis of the operators~$H'_{\al}$, $\al \in S^1$. The
eigenvalue of $H'_{\al}$ on $f_n$ is ${\rm e}^{2\pi{\rm i} n\al}$.

Likewise, in the two-dimensional case of the elliptic curve $E_{\rm i}$, the
Fourier harmonics $f_{m,n}$ form an orthogonal basis of eigenfunctions
of the operators $_{p_0}H_p$, $p \in E_{\rm i}$, in $L^2(E_{\rm i})$:
\begin{gather*} 
_{p_0}H_p \cdot f_{m,n} = {\rm e}^{2\pi{\rm i} (mx_p+ny_p)} f_{m,n}.
\end{gather*}
From this formula we see that the eigenvalue of $_{p_0}H_p$ on
$f_{m,n}$ is ${\rm e}^{2\pi{\rm i} (mx_p+ny_p)}$. Thus, we have obtained a complete
description of the Hecke eigenfunctions and eigenvalues for the curve
$X=E_{\rm i}$ and the group $G={\rm GL}_1$.

\subsection{General elliptic curve} \label{genell}

Now we generalize this to the case of an arbitrary elliptic curve
$E_\tau \simeq \C/(\Z+\Z\tau)$ with $\on{Im}\tau>0$. Recall that we
identify every component of $\Pic(E_\tau)$ with $E_\tau$ using the
reference point $p_0$. Then we obtain the Hecke operators $_{p_0}H_p$
labeled by $p \in E_\tau$ given by the shift by $p$ naturally acting
on $E_\tau$ (see formula \eqref{Hpdef}). The eigenfunctions and
eigenvalues of these operators are then given by the following
theorem.

\begin{thm} \label{Etau eig}
The joint eigenfunctions of the Hecke operators $_{p_0}H_p$, $p \in
E_\tau$, on $L^2(E_\tau)$ are
\begin{gather*} 
f^\tau_{m,n}(z,\ol{z}) = {\rm e}^{2\pi
 {\rm i}m(z\ol\tau-\ol{z}\tau)/(\ol\tau-\tau)} \cdot {\rm e}^{2\pi
 {\rm i}n (z-\ol{z})/(\tau-\ol\tau)}, \qquad m,n \in \Z.
\end{gather*}
The eigenvalues are given by the right hand side of the following
formula:
\begin{gather} \label{eigmntau}
_{p_0}H_p \cdot f^\tau_{m,n} = {\rm e}^{2\pi
 {\rm i}m(p\ol\tau-\ol{p}\tau)/(\ol\tau-\tau)} \cdot {\rm e}^{2\pi
 {\rm i}n (p-\ol{p})/(\tau-\ol\tau)} f^\tau_{m,n}.
\end{gather}
\end{thm}

In Section~\ref{higher} we will give an alternative formula for these eigenfunctions (for an arbitrary smooth projective curve instead of~$E_\tau$).

\subsection[Digression: Eigenvalues of the Hecke operators and representations of the fundamental group]{Digression: Eigenvalues of the Hecke operators\\ and representations of the fundamental group} \label{fund}

Let ${\mb H}(E_\tau)$ be the spectrum of the algebra of Hecke
operators acting on $L^2\big(\Pic^0(E_\tau)\big) = L^2(E_\tau)$. In this
subsection we compare the description of ${\mb H}(E_\tau)$ given in
Theorem~\ref{Etau eig} with that envisioned by Langlands in~\cite{L:analyt}.

Let ${\mb E}(E_\tau)$ be the set of equivalence classes of
one-dimensional representations of the fundamental group
$\pi_1(E_\tau,p_0)$ with finite image. In~\cite{L:analyt}, Langlands
attempts to construct a one-to-one correspondence between ${\mb
 H}(E_\tau)$ and ${\mb E}(E_\tau)$ in two different ways.

The first is to express the Hecke eigenvalues corresponding to a given
Hecke eigenfunction as holonomies of a flat unitary connection on a
line bundle on $E_\tau$ and then take the monodromy representation of
this connection (see part~(4) in Section~1.4). I show below
that it is indeed possible to express the Hecke eigenvalues that we
have found in Theorem~\ref{Etau eig} as holonomies of a flat unitary
connection on the trivial line bundle on $E_\tau$ (furthermore, this
will be generalized in Section~\ref{higher} to the case of an
arbitrary curve~$X$). But all of these connections have {\em trivial}
monodromy representation. Thus, the map ${\mb H}(E_\tau) \to {\mb
 E}(E_\tau)$ we obtain this way is trivial, i.e., its image consists
of a single element of ${\mb E}(E_\tau)$. (Herein lies an important
difference between the analytic and geometric theories for curves over~$\C$, which is discussed in more detail in Remark~\ref{rem1} below.)

Second, Langlands attempted to construct a map ${\mb H}(E_\tau) \to
{\mb E}(E_\tau)$ explicitly. Unfortunately, this construction does not
yield a bijective map, either, as I show in Remark~\ref{biject} below.

Let me show how to express the eigenvalues of the Hecke operators
$_{p_0}H_p$, $p \in E_\tau$, on a given eigenfunction as holonomies of a
flat unitary connection.

Consider first the case of $\tau={\rm i}$. In this case, we assign to the
Hecke eigenfunction $f_{m,n}$ given by formula~\eqref{fmn} the
following unitary flat connection~$\nabla^{(m,n)}$ on the trivial line
bundle over $E_{\rm i}$:
\begin{gather*} \label{nablamn}
\nabla^{(m,n)} = d - 2\pi{\rm i} m \, {\rm d}x - 2\pi{\rm i} n \, {\rm d}y
\end{gather*}
(since the line bundle is trivial, a connection on it is the same as a
one-form on the curve). In other words, the corresponding first order
differential operators along~$x$ and~$y$ are given by the formulas
\begin{gather*} 
\nabla^{(m,n)}_x = \frac{\pa}{\pa x} - 2\pi{\rm i} m,\qquad
\nabla^{(m,n)}_y = \frac{\pa}{\pa y} - 2\pi{\rm i} n.
\end{gather*}
The horizontal sections of this connection are the solutions of the
equations
\begin{gather} \label{diffeq}
\nabla^{(m,n)}_x \cdot \Phi = \nabla^{(m,n)}_y \cdot \Phi = 0.
\end{gather}
They have the form
\begin{gather*}
\Phi_{m,n}(x,y) = {\rm e}^{2\pi{\rm i} (mx+ny)}
\end{gather*}
up to a scalar. The function $\Phi_{m,n}$ is the unique solution of~\eqref{diffeq} normalized so that its value at the point $0 \in E_{\rm i}$,
corresponding to our reference point $p_0 \in E_{\rm i}$, is equal to~$1$.
The value of this function $\Phi_{m,n}$ at $p=x_p+{\rm i}y_p \in \C/(\Z+\Z
{\rm i})$ is indeed equal to the eigenvalue of the Hecke operator
$_{p_0}H_p$ on the harmonic $f_{m,n}$.

Thus, this eigenvalue can be represented as the holonomy of the
connection~$\nabla^{(m,n)}$ over a~path connecting our reference point
$p_0 \in E_{\rm i}$, which corresponds to $0 \in \C/(\Z+\Z{\rm i})$, and the
point $p \in E_{\rm i}$. Since the connection is flat, it does not matter
which path we choose.

However, and this is a crucial point, the connection
$\nabla^{(m,n)}$ has {\em trivial monodromy} on $E_{\rm i}$. Indeed,
\begin{gather*}
\Phi_{m,n}(x+1,y) = \Phi_{m,n}(x,y+1) = \Phi_{m,n}(x,y)
\end{gather*}
for all $m,n \in \Z$.

Similarly, we assign a flat unitary connection $_\tau\nabla^{(m,n)}$
on the trivial line bundle on $E_\tau$ for each Hecke eigenfunction
$f^\tau_{m,n}$:
\begin{gather*} 
_\tau\nabla^{(m,n)} = d - 2\pi{\rm i}
\frac{n-m\ol\tau}{\tau-\ol\tau} \,{\rm d}z - 2 \pi{\rm i}
\frac{m\tau-n}{\tau-\ol\tau} \,{\rm d}\ol{z}.
\end{gather*}
The first order operators corresponding to $z$ and $\ol{z}$ are
\begin{gather*} 
_\tau\nabla^{(m,n)}_z = \frac{\pa}{\pa z} - 2\pi{\rm i}\frac{n-m\ol\tau}{\tau-\ol\tau},\\ 
_\tau\nabla^{(m,n)}_{\ol{z}} = \frac{\pa}{\pa \ol{z}} - 2 \pi{\rm i} \frac{m\tau-n}{\tau-\ol\tau}.
\end{gather*}
Just as in the case $\tau={\rm i}$, for every $p \in E_\tau$, the holonomy
of the connection $_\tau\nabla^{(m,n)}$ over a path connecting $p_0
\in E_\tau$ and $p \in E_\tau$ is equal to the eigenvalue of $_{p_0}
H_p$ on $f^\tau_{m,n}$ given by the right hand side of formula
\eqref{eigmntau}. However, as in the case of $\tau={\rm i}$, all connections
$_\tau\nabla^{(m,n)}$ yield the trivial monodromy representation
$\pi_1(E_\tau,p_0) \to {\rm GL}_1$.

\begin{rem} \label{biject}
On pp.~59--60 of~\cite{L:analyt}, another attempt is made to construct
a map from the set~${\mb H}(E_\tau)$ (the spectrum of the algebra of Hecke
operators acting on~$L^2(E_\tau)$) to the set ${\mb E}(E_\tau)$ of
equivalence classes of homomorphisms $\pi_1(E_\tau,p_0) \to {\rm GL}_1$ with
finite image. According to Theorem~\ref{Etau eig}, the set ${\mb
 H}(E_\tau)$ is identified with $\Z \times \Z$. On the other
hand, the set ${\mb E}(E_\tau)$ can be identified with $\mu
\times \mu$, where $\mu$ is the group of complex
roots of unity (we have an isomorphism $\Q/\Z \simeq \mu$ sending $\ka
\in \Q/\Z$ to ${\rm e}^{2\pi{\rm i} \kappa}$). Indeed, since $\pi_1(E_\tau,p_0)
\simeq \Z \times \Z$, a homomorphism
$\phi\colon \pi_1(E_\tau,p_0) \to {\rm GL}_1 \simeq \C^\times$ is uniquely
determined by its values on the elements $A = (1,0)$ and $B = (0,1)$
of $\Z \times \Z$. The homomorphism $\phi$ has finite image if and
only if both $\phi(A)$, $\phi(B)$ belong to~$\mu$.

Langlands attempts to construct a map $(\Z \times \Z) \to (\mu \times
\mu)$ as follows (see pp.~59--60 of~\cite{L:analyt}): he sets
\begin{gather*} 
(0,0) \mapsto (1,1).
\end{gather*}
Next, given a non-zero element $(k,l) \in \Z \times \Z$, there exists
a matrix $g_{k,l} = \begin{pmatrix} \alpha & \beta \\ \gamma &
 \delta \end{pmatrix} \in {\rm SL}_2(\Z)$ such that
\begin{gather} \label{matrix}
\begin{pmatrix} k & l \end{pmatrix} = \begin{pmatrix} k' &
 0 \end{pmatrix} \begin{pmatrix} \alpha & \beta \\
 \gamma & \delta \end{pmatrix}, \qquad k'>0.
\end{gather}
Two comments on \eqref{matrix}: first, as noted in~\cite{L:analyt},
the matrix $g_{k,l}$ is not uniquely determined by formula~\eqref{matrix}. Indeed, this formula will still be satisfied if we
multiply $g_{k,l}$ on the left by any lower triangular matrix in
${\rm SL}_2(\Z)$. Second, formula~\eqref{matrix} implies that
\begin{gather} \label{map1}
(k,l) = k'(\alpha,\beta), \qquad \on{gcd}(\alpha,\beta)=\pm 1, \qquad
k'>0,
\end{gather}
where, for a pair of integers $(k,l) \neq (0,0)$, we define
$\on{gcd}(k,l)$ as $l$ if $k=0$, as $k$ if $l=0$, and
$\on{gcd}(|k|,|l|)$ times the product of the signs of $k$ and $l$ if
they are both non-zero. Therefore
\begin{gather*} 
k' = |\on{gcd}(k,l)|.
\end{gather*}

Using a particular choice of the matrix $g_{k,l}$, Langlands defines a
new set of generators $\{ A',B' \}$ of the group $\pi_1(E_\tau,p_0)$:
\begin{gather} \label{change}
A' = A^\alpha B^\beta, \qquad B' = A^\ga B^\delta.
\end{gather}
He then defines a homomorphism $\phi_{k,l}\colon \pi_1(E_\tau,p_0) \to {\rm GL}_1$ corresponding to $(k,l)$ by the formulas
\begin{gather*}
A' \mapsto {\rm e}^{2\pi{\rm i}/k'}, \qquad B' \mapsto 1.
\end{gather*}
Now, formula \eqref{change} implies that
\begin{gather*} 
A = (A')^\delta (B')^{-\beta}, \qquad B = (A')^{-\ga} (B')^\al,
\end{gather*}
and so we find the values of $\phi_{k,l}$ on the original generators~$A$ and~$B$:
\begin{gather} \label{rho}
A \mapsto {\rm e}^{2\pi{\rm i} \delta/k'}, \qquad B \mapsto {\rm e}^{-2\pi{\rm i} \ga/k'}.
\end{gather}
Langlands writes in~\cite{L:analyt}, ``This has a peculiar property that part of the numerator becomes the denominator, which baffles me and may well baffle the reader''. He goes on to say, ``To be honest,
this worries me''.

In fact, this construction does {\em not} give us a well-defined map
$(\Z \times \Z) \to (\mu \times \mu)$. Indeed, $g_{k,l}$ is only
defined up to left multiplication by a lower triangular matrix:
\begin{gather*} 
\begin{pmatrix} \alpha & \beta \\
 \gamma & \delta \end{pmatrix} \mapsto \begin{pmatrix} 1
 & 0 \\ x & 1 \end{pmatrix} \begin{pmatrix}
 \alpha & \beta \\ \gamma & \delta \end{pmatrix}, \qquad x \in \Z,
\end{gather*}
under which we have the following transformation:
\begin{gather*} 
\gamma \mapsto \gamma+x\alpha, \qquad \delta \mapsto \delta + x\beta.
\end{gather*}
But then the homomorphism \eqref{rho} gets transformed to the
homomorphism sending
\begin{gather} \label{rho1}
A \mapsto {\rm e}^{2\pi{\rm i} (\delta+x\beta)/k'}, \qquad B \mapsto {\rm e}^{-2\pi{\rm i}
 (\ga+x\alpha)/k'}.
\end{gather}
The homomorphisms \eqref{rho} and \eqref{rho1} can only coincide for
all $x \in \Z$ if both $\alpha$ and $\beta$ are divisible by $k'$. But
if $k' \neq 1$, this contradicts the condition, established in formula
\eqref{map1}, that $\alpha$ and $\beta$ are relatively prime. Hence
\eqref{rho} and \eqref{rho1} will in general differ from each other,
and so we don't get a well-defined map $(\Z \times \Z) \to (\mu \times
\mu)$.

We could try to fix this problem by replacing the relation~\eqref{change} with
\begin{gather*} 
A = (A')^\alpha (B')^\gamma, \qquad B = (A')^\beta (B')^\delta.
\end{gather*}
Then the homomorphism $\phi_{k,l}$ would send
\begin{gather*} 
A \mapsto {\rm e}^{2\pi{\rm i} \alpha/k'}, \qquad B \mapsto {\rm e}^{2\pi{\rm i} \beta/k'}.
\end{gather*}
This way, we get a well-defined map $(\Z \times \Z) \to (\mu \times
\mu)$, but it's not a bijection.
\end{rem}

In fact, there is no reason to expect that there is a meaningful
bijection between the above sets ${\mb H}(E_\tau)$ and ${\mb
 E}(E_\tau)$. Indeed, according to Theorem~\ref{Etau eig}, the set
${\mb H}(E_\tau)$ can be naturally identified with the group of
continuous characters $E_\tau \to \C^\times$ (where $E_\tau$ is viewed
as an abelian group), which is isomorphic to $\Z \times \Z$.

On the other hand, let ${\mb E}(E_\tau)$ is the subgroup of elements
of finite order in the group of cha\-rac\-ters $\pi_1(E_\tau,p_0) \to
\C^\times$. The whole group of such characters, which is isomorphic to
\mbox{$\C^\times \times \C^\times$}, is the {\em dual group} of $\Z \times \Z
= {\mb H}(E_\tau)$. The set~${\mb E}(E_\tau)$ is its subgroup of
elements of finite order, which isomorphic to~$\mu \times \mu$, where~$\mu$ is the (multiplicative) group of complex roots of unity. Clearly, $\Z \times \Z$ and $\mu \times \mu$ are not isomorphic
as abstract groups. Of course, since each of these two sets is
countable, there exist bijections between them as sets. But it's hard
to imagine that such a bijection would be pertinent to the questions
at hand.

\begin{rem} \label{rem1} Recall that in the classical unramified Langlands correspondence
for a curve over~$\Fq$, to each joint eigenfunction of the Hecke
operators we assign a Langlands parameter. In the case of
$G={\rm GL}_n$, this is an equivalence class of $\ell$-adic homomorphisms
from the \'etale fundamental group of~$X$ to~${\rm GL}_n$
(and more generally, one considers homomorphisms to the Langlands dual
group~$\LG$ of~$G$). Given such a homomorphism $\sigma$, to each closed
point~$x$ of~$X$ we can assign an $\ell$-adic number, the trace of
$\sigma(\on{Fr}_x)$, where~$\on{Fr}_x$ is the Frobenius conjugacy
class, so we obtain a function from the set of closed points of $X$ to
the set of conjugacy classes in~${\rm GL}_n(\ol{\mathbb Q}_\ell)$.

In the geometric Langlands correspondence for curves over $\C$, the
picture is different. Now the role of the \'etale fundamental group is
played by the topological fundamental group~$\pi_1(X,p_0)$. Thus, the
Langlands parameters are the equivalence classes of homomorphisms
$\pi_1(X,p_0) \to {\rm GL}_n$ (or, more generally, to~$\LG$). The question
then is: how to interpret such a~homomorphism as a~Hecke
``eigenvalue'' on a~Hecke eigensheaf?

The point is that for a Hecke eigensheaf, the ``eigenvalue'' of a
Hecke operator (or rather, Hecke functor) is not a number but an
$n$-dimensional {\em vector space}. As we move along a closed path on
our curve (starting and ending at the point $p_0$ say), this vector
space will in general undergo a non-trivial linear transformation,
thus giving rise to a non-trivial homomorphism $\pi_1(X,p_0) \to
{\rm GL}_n$.

Note that over $\C$ we have the Riemann--Hilbert correspondence, which
sets up a bijection between the set of equivalence classes of
homomorphisms $\pi_1(X,p_0) \to {\rm GL}_n$ (or, more generally,
$\pi_1(X,p_0) \to \LG$) and the set of equivalence classes of pairs
$({\mc P},\nabla)$, where ${\mc P}$ is a rank $n$ bundle on $X$ (or,
more generally, an $\LG$-bundle) and $\nabla$ is a flat connection on
${\mc P}$. The map between the two data is defined by assigning to
$({\mc P},\nabla)$ the monodromy representation of $\nabla$
(corresponding to a specific a trivialization of ${\mc P}$ at
$p_0$). We may therefore take equivalence classes of the flat bundles
$({\mc P},\nabla)$ as our Langlands parameters instead of equivalence
classes of homomorphisms $\pi_1(X,p_0) \to {\rm GL}_n$. As explained in the
previous paragraph, these flat bundles $({\mc P},\nabla)$ will in
general have non-trivial monodromy.

However, in this section we consider (in the case of ${\rm GL}_1$ and a
curve $X$) the {\em eigenfunctions} of the Hecke operators $_{p_0}H_p$,
$p \in X$, on $\Pic^0(X)$. Their eigenvalues are {\em numbers}, not
vector spaces. Therefore they cannot undergo any transformations as we
move along a closed path on our curve. In other words, these numbers
give rise to a {\em single-valued} function from $X$ to ${\rm GL}_1(\C)$ (it
actually takes values in $U_1 \subset {\rm GL}_1(\C)$). Because the function
is single-valued, if we represent this function as the holonomy of a
flat connection on a line bundle on $X$, then this connection
necessarily has trivial monodromy. And indeed, we have seen above that
each collection of joint eigenvalues of the Hecke operators
$_{p_0}H_p$, $p \in E_\tau$, on functions on~$\Pic^0(E_\tau)$ can be
represented as holonomies of a specific (unitary) connection
$_\tau\nabla^{m,n}$ with trivial monodromy. The same is true for other
curves, as we will see below.
\end{rem}

\subsection{Higher genus curves} \label{higher}

Let $X$ be a smooth projective connected curve over $\C$. Denote by
$\Pic(X)$ the Picard variety of~$X$, i.e., the moduli space of line
bundles on~$X$ (as before, the moduli stack $\Bun_{{\rm GL}_1}(X)$ of line
bundles on~$X$ is the quotient of $\Pic(X)$ by the trivial action of
${\mathbb G}_m={\rm GL}_1$). We have a~decomposition of $\Pic(X)$ into a~disjoint union of connected components $\Pic^d(X)$ corresponding to line bundles of degree~$d$. The Hecke operator~$H_p$, $p \in X$, is the pull-back of functions with respect to the map (see formula~\eqref{origT} for $X=E_\tau$):
\begin{align*} 
T_p\colon \ \Pic^d(X) &\to \Pic^{d+1}(X), \\
\notag {\mc L} &\mapsto {\mc L}(p).
\end{align*}

The Hecke operators $H_p$ with different $p \in X$ commute with each
other, and it is natural to consider the problem of finding joint
eigenfunctions and eigenvalues of these operators on functions on
$\Pic(X)$. In the same way as in Section \ref{anell}, we find that
this problem is equivalent to the problem of finding joint
eigenfunctions and eigenvalues of the operators $_{p_0}H_p =
H_{p_0}^{-1} H_p$ on functions on $\Pic^0(X)$, where $p_0$ is a~reference point on~$X$ that we choose once and for all. The operator
$_{p_0}H_p$ is the pull-back of functions with respect to the map
$_{p_0}T_p\colon \Pic^0(X) \to \Pic^0(X)$ sending a line bundle ${\mc L}$
to ${\mc L}(p-p_0)$.

Now, $\Pic^0(X)$ is the Jacobian of $X$, which is a $2g$-dimensional
torus (see, e.g., \cite{GH})
\begin{gather*}
\Pic^0(X) \simeq H^0\big(X,\Omega^{1,0}\big)^*/H_1(X,\Z),
\end{gather*}
where $H_1(X,\Z)$ is embedded into the space of linear functionals on the space $H^0\big(X,\Omega^{1,0}\big)$ of holomorphic one-forms on~$X$ by sending $\beta \in H_1(X,\Z)$ to the linear functional
\begin{gather} \label{embed hol}
\omega \in H^0\big(X,\Omega^{1,0}\big) \mapsto \int_\beta \omega.
\end{gather}

Motivated by Theorem \ref{Etau eig}, it is natural to guess that the
standard Fourier harmonics in $L^2\big(\Pic^0(X)\big)$ form an orthogonal
eigenbasis of the Hecke operators. This is indeed the case.

To see that, we give an explicit formula for these harmonics. They can
be written in the form ${\rm e}^{2 \pi{\rm i} \varphi}$, where $\varphi\colon
H^0\big(X,\Omega^{1,0}\big)^* \to \R$ is an $\R$-linear functional such that
$\varphi(\beta) \in \Z$ for all $\beta \in H_1(X,\Z)$. To write them
down explicitly, we use the Hodge decomposition
\begin{gather*}
H^1(X,\C) = H^0\big(X,\Omega^{1,0}\big) \oplus H^0\big(X,\Omega^{0,1}\big) =
H^0\big(X,\Omega^{1,0}\big) \oplus \ol{H^0\big(X,\Omega^{1,0}\big)}
\end{gather*}
to identify $H^0\big(X,\Omega^{1,0}\big)$, viewed as an $\R$-vector
space, with $H^1(X,\R)$ by the formula
\begin{gather} \label{ident vsp}
\omega \in H^0\big(X,\Omega^{1,0}\big) \mapsto \omega + \ol\omega.
\end{gather}
In particular, for any class $c \in H^1(X,\R)$, there is a unique
holomorphic one-form $\omega_c$ such that~$c$ is represented by the
real-valued harmonic one-form $\omega_c + \ol{\omega}_c$,
\begin{gather} \label{omegac}
H^1(X,\R) \ni c = \omega_c + \ol{\omega}_c ,
\qquad \omega_c \in H^0\big(X,\Omega^{1,0}\big).
\end{gather}

Viewed as a real manifold,
\begin{gather*}
\Pic^0(X) \simeq H^1(X,\R)^*/H_1(X,\Z),
\end{gather*}
where $H_1(X,\Z)$ is embedded into $H^1(X,\R)^*$ by sending $\beta \in
H_1(X,\Z)$ to the linear functional on~$H^1(X,\R)$ given by the
formula (compare with formulas~\eqref{embed hol} and~\eqref{ident vsp})
\begin{gather*} 
H^1(X,\R) \ni c \mapsto \int_\beta c = \int_\beta (\omega_c + \ol{\omega}_c).
\end{gather*}

Now, to each $\gamma \in H^1(X,\Z)$ we attach the corresponding
element of the vector space $H^1(X,\R)$, which can be viewed as a
linear functional $\varphi_\ga$ on the dual vector space
$H^1(X,\R)^*$,
\begin{gather*}
\varphi_\ga\colon \ H^1(X,\R)^* \to \R.
\end{gather*}
It has the desired property: $\varphi_\ga(\beta) \in \Z$ for all
$\beta \in H_1(X,\Z)$. The corresponding functions
\begin{gather} \label{harmonics}
{\rm e}^{2\pi{\rm i} \varphi_\ga}, \qquad \gamma \in H^1(X,\Z),
\end{gather}
are the Fourier harmonics that form an orthogonal basis of the Hilbert space $L^2\big(\Pic^0(X)\big)$.

We claim that each of these functions is an eigenfunction of the Hecke
operators $_{p_0}H_p$, $p \in X$, so that together they give us a
sought-after orthogonal eigenbasis of the Hecke operators. To see that, we use the Abel--Jacobi map.

For $d>0$, let $X^{(d)}$ be the $d$th symmetric power of $X$, and
$p_d\colon X^{(d)} \to \Pic^d(X)$ the Abel--Jacobi map
\begin{gather*} 
p_d(D) = \OO(D), \qquad D = \sum_{i=1}^d [x_i], \qquad x_i \in X.
\end{gather*}
We can lift the map $T_p$ to a map
\begin{align*} 
\wt{T}_p\colon \ X^{(d)} &\to X^{(d+1)}, \\ \notag
D &\mapsto D+[p],
\end{align*}
so that we have a commutative diagram
\begin{gather} \label{compat}
\begin{CD}
X^{(d)} @>{\wt{T}_p}>> X^{(d+1)} \\
@Vp_dVV @VVp_{d+1}V \\
\Pic^d(X) @>T_p>> \Pic^{d+1}(X).
\end{CD}
\end{gather}
Denote by $\wt{H}_p$ the corresponding pull-back operator on functions.

Now let $f_0$ be a non-zero function on $\Pic^0(X)$. Identifying
$\Pic^d(X)$ with $\Pic^0(X)$ using the reference point $p_0$:
\begin{gather} \label{Picd0}
{\mc L} \mapsto {\mc L}(-d \cdot p_0),
\end{gather}
we obtain a non-zero function $f_d$ on $\Pic^d(X)$ for all $d \in
\Z$. Let $\wt{f}_d$ the pull-back of $f_d$ to $X^{(d)}$ for
$d>0$. Suppose that these functions satisfy
\begin{gather*} 
\wt{H}_p\big(\wt{f}_{d+1}\big) = \lambda_p \wt{f}_d, \qquad
 p \in X, \qquad d>0,
\end{gather*}
where $\la_p \neq 0$ for all $p$ and $\la_{p_0}=1$. This is equivalent
to the following factorization formula for~$\wt{f}_d$:
\begin{gather} \label{wtfd}
\wt{f}_d \left( \sum_{i=1}^d [x_i] \right) = c \prod_{i=1}^d
\la_{x_i}, \qquad c \in \C, \qquad d>0.
\end{gather}
The surjectivity of $p_d$ with $d\geq g$ and the commutativity of the
diagram \eqref{compat} then implies that
\begin{gather*} 
H_p(f_{d+1}) = \lambda_p f_d, \qquad p \in X, \qquad d\geq g.
\end{gather*}
But then it follows from the definition of $f_d$ that $f_0$ is an eigenfunction of the operators $_{p_0}H_p = H_{p_0}^{-1} H_p$ with the eigenvalues $\la_p = \wt{f}_1([p])$.

This observation gives us an effective way to demonstrate that a given function $f_0$ on $\Pic^0(X)$ is a Hecke eigenfunction.

Let us use it in the case of the function $f_0={\rm e}^{2\pi{\rm i}
 \varphi_{\ga}}$, $\ga \in H^1(X,\Z)$, on $\Pic^0(X)$ given by formula
\eqref{harmonics}. For that, denote by
\begin{gather} \label{d exp}
_d {\rm e}^{2\pi{\rm i} \varphi_{\ga}}, \qquad \ga \in H^1(X,\Z),
\end{gather}
the corresponding functions $f_d$ on $\Pic^d(X)$ obtained via the
identification \eqref{Picd0}. We claim that for any $\ga \in
H^1(X,\Z)$, the pull-backs of $_d {\rm e}^{2\pi{\rm i} \varphi_{\ga}}$ to
$X^{(d)}$, $d>0$, via the Abel--Jacobi maps have the form \eqref{wtfd},
and hence ${\rm e}^{2\pi{\rm i} \varphi_{\ga}}$ is a Hecke eigenfunction on
$\Pic^0(X)$.

To see that, we recall an explicit formula for the composition
\begin{gather} \label{XdOmega}
X^{(d)} \to \Pic^d(X) \to \Pic^0(X) \simeq
H^0\big(X,\Omega^{1,0}\big)^*/H_1(X,\Z),
\end{gather}
where the second map is given by formula \eqref{Picd0} (see, e.g.,
\cite{GH}). Namely, the composition \eqref{XdOmega} maps $\sum\limits_{i=1}^d
[x_i] \in X^{(d)}$ to the linear functional on $H^0\big(X,\Omega^{1,0}\big)$
sending
\begin{gather*}
\omega \in H^0\big(X,\Omega^{1,0}\big) \mapsto \sum_{i=1}^d \int_{p_0}^{x_i} \omega.
\end{gather*}
Composing the map~\eqref{XdOmega} with the isomorphism
$H^0\big(X,\Omega^{1,0}\big) \simeq H^1(X,\R)$ defined above, we obtain a map
\begin{gather*} 
_{p_0} \Phi_d\colon \ X^{(d)} \to H^1(X,\R)^*/H_1(X,\Z),
\end{gather*}
which maps $\sum\limits_{i=1}^d [x_i] \in X^{(d)}$ to the linear functional
$_{p_0}\Phi_d \Big(\sum\limits_{i=1}^d [x_i]\Big)$ on $H^1(X,\R)$ given by
the formula
\begin{gather*} 
_{p_0}\Phi_d\left(\sum_{i=1}^d [x_i]\right)\colon \ c \in H^1(X,\R) \mapsto
\sum_{i=1}^d \int_{p_0}^{x_i} (\omega_c + \ol{\omega}_c)
\end{gather*}
(see formula \eqref{omegac} for the definition of $\omega_c$).

If $c \in H^1(X,\R)$ is the image of an {\em integral}
cohomology class
\begin{gather*}
\ga \in H^1(X,\Z),
\end{gather*}
we will write the corresponding holomorphic one-form $\omega_c$ as $\omega_\ga$.

Let $_{p_0}\wt{f}_{d,\ga}$ be the pull-back of the function $_d
{\rm e}^{2\pi{\rm i} \varphi_{\ga}}$ (see formula~\eqref{d exp}) to~$X^{(d)}$. Equivalently, $_{p_0}\wt{f}_{d,\ga}$ is the pull-back of
the function ${\rm e}^{2\pi{\rm i} \varphi_{\ga}}$ under the map
$_{p_0}\Phi_d$. It follows from the definition of $_{p_0}\Phi_d$ that
the value of $_{p_0}\wt{f}_{d,\ga}$ at $\sum\limits_{i=1}^d [x_i]$ is equal
to
\begin{gather*}
\exp\left( 2\pi{\rm i} \; \; {}_{p_0}\Phi_d\left(\sum_{i=1}^d
 [x_i]\right)(\ga) \right) = \exp\left( 2\pi{\rm i} \; \; \sum_{i=1}^d
 \int_{p_0}^{x_i} (\omega_\ga + \ol{\omega}_\ga) \right).
\end{gather*}

Thus, we obtain that $_{p_0}\wt{f}_{d,\ga}$ is given by the formula
\begin{gather} \label{p0 f}
_{p_0}\wt{f}_{d,\ga} \left( \sum_{i=1}^d [x_i] \right) = \exp\left( 2
 \pi{\rm i} \; \; \sum_{i=1}^d \int_{p_0}^p (\omega_\ga + \ol{\omega}_\ga)
\right) = \prod_{i=1}^d \la^\ga_{x_i},
\end{gather}
where
\begin{gather} \label{laga}
\la^\ga_p = {\rm e}^{2 \pi{\rm i} \int_{p_0}^p (\omega_\ga + \ol{\omega}_\ga)}.
\end{gather}

We conclude that the functions $_{p_0}\wt{f}_{d,\ga}$ satisfy the
factorization property~\eqref{wtfd}. Therefore the function ${\rm e}^{2\pi{\rm i}
 \varphi_{\ga}}$ on $\Pic^0(X)$ is indeed an eigenfunction of~$_{p_0}H_p$, with the eigenvalue $\la^\ga_p$ given by formula~\eqref{laga}, which is what we wanted to prove.\footnote{Note that
 Abel's theorem implies that each function $_{p_0}\wt{f}_{d,\ga}$, $\ga
 \in H^1(X,\Z)$, is constant along the fibers of the Abel--Jacobi map
 $X^{(d)} \to \Pic^d$ and therefore descends to $\Pic^d$. This
 suggests another proof of Theorem~\ref{X eig}: we start from the
 functions $_{p_0}\wt{f}_{d,\ga}$ on $X^{(d)}$, $d>0$. Formula~\eqref{p0 f} shows that they combine into an eigenfunction of the
 operators $\wt{H}_p$. Hence the function on $\Pic^d(X)$, $d\geq g$, to
 which $_{p_0}\wt{f}_{d,\ga}$ descends, viewed as a function on~$\Pic^0(X)$ under the identification~\eqref{Picd0}, is a Hecke
 eigenfunction. One can then show that this function is equal to~${\rm e}^{2\pi{\rm i} \varphi_{\ga}}$.}

Thus, we have proved the following theorem.

\begin{thm} \label{X eig} The joint eigenfunctions of the Hecke operators $_{p_0}H_p$, $p \in
X$, on $L^2\big(\Pic^0(X)\big)$ are the functions ${\rm e}^{2\pi{\rm i} \varphi_{\ga}}$,
$\ga \in H^1(X,\Z)$. The eigenvalues of~$_{p_0}H_p$ are given by
formula \eqref{laga}, so that we have
\begin{gather*} 
_{p_0}H_p \cdot {\rm e}^{2\pi{\rm i} \varphi_{\ga}} = {\rm e}^{2 \pi{\rm i} \int_{p_0}^p
 (\omega_\ga + \ol{\omega}_\ga)} {\rm e}^{2\pi{\rm i} \varphi_{\ga}}.
\end{gather*}
\end{thm}

As in the case of an elliptic curve discussed in Section~\ref{genell},
the eigenvalues~\eqref{laga} can be interpreted as the holonomies of
the flat unitary connections
\begin{gather*}
\nabla_\ga = d - 2\pi{\rm i}(\omega_\ga + \ol\omega_\ga), \qquad \ga \in H^1(X,\Z)
\end{gather*}
on the trivial line bundle on~$X$, taken along (no matter which) path
from~$p_0$ to~$p$. As in the case of elliptic curves, the monodromy
representation of each of these connections is {\em trivial}, ensuring
that the Hecke eigenvalues $\la^\ga_p$, viewed as functions of $p \in
X$, are single-valued (see Section~\ref{fund}).

\subsection{General torus} \label{torus}

Let now $T$ be a connected torus over $\C$, and $\Bun_T(X)$ the moduli
space of $T$-bundles on~$X$ (note that the moduli stack $\Bun_T(X)$ is
the quotient of $\Bun_T(X)$ by the trivial action of~$T$). In Section~\ref{higher} we find the joint eigenfunctions and eigenvalues of the Hecke operators in the case of $\Bun_T(X)$ where $T={\mathbb G}_m$; in
this case $\Bun_{{\mathbb G}_m}(X) = \Pic(X)$. Here we generalize these
results to the case of an arbitrary~$T$.

Let $\Lambda^*(T)$ and $\Lambda_*(T)$ be the lattices of characters
and cocharacters of $T$, respectively. Any ${\mc P} \in \Bun_T(X)$ is
uniquely determined by the ${\mathbb G}_m$-bundles (equivalently, line
bundles) ${\mc P} \underset{{\mathbb G}_m}\times \chi$ associated to
the characters $\chi\colon T \to {\mathbb G}_m$ in $\Lambda^*(T)$. This
yields a canonical isomorphism
\begin{gather*}
\Bun_T(X) \simeq \Pic(X) \underset{\Z}\otimes \Lambda_*(T) =
\bigsqcup_{\check\nu \in \Lambda_*(T)} \Bun^{\check\nu}_T(X).
\end{gather*}
The neutral component
\begin{gather*}
\Bun^0_T(X) = \Pic^0(X) \underset{\Z}\otimes \Lambda_*(T)
\end{gather*}
is non-canonically isomorphic to $\Pic^0(X)^r$, where $r$ is the rank
of the lattice $\Lambda_*(T)$.

The Hecke operators $H_p^{\check\mu}$ are now labeled by $p \in X$ and
$\check\mu \in \Lambda_*(T)$. The operator $H_p^{\check\mu}$ corresponds to
the pull-back under the map
\begin{align*} 
T_p^{\check\mu}\colon \ \Bun^{\check\nu}_T(X) &\to
\Bun^{\check\nu+\check\mu}_T(X), \\ \notag
{\mc P} &\mapsto {\mc P}(\check\mu \cdot p),
\end{align*}
where ${\mc P}(\check\mu \cdot p)$ is defined by the formula
\begin{gather*}
{\mc P}(\check\mu \cdot p) \underset{{\mathbb G}_m}\times \chi =
({\mc P} \underset{{\mathbb G}_m}\times \chi)(\langle \chi,\check\mu
\rangle \cdot p), \qquad \chi \in \Lambda^*(T).
\end{gather*}
As in the case of $T={\mathbb G}_m$, we choose, once and for all, a
reference point $p_0 \in X$.

As in the case of $T={\mathbb G}_m$, finding eigenfunctions and
eigenvalues of the commuting operators~$H_p^{\check\mu}$ on functions
on $\Bun_T(X)$ is equivalent to finding eigenfunctions and eigenvalues
of the operators
\begin{gather*}
_{p_0}H_p^{\check\mu} = \big(H_{p_0}^{\check\mu}\big)^{-1} \circ
H_p^{\check\mu}
\end{gather*}
on functions on $\Bun^0_T(X)$. As in the case of $T={\mathbb G}_m$, we
represent $\Bun^0_T(X)$ as
\begin{gather} \label{JacT}
\Bun^0_T(X) \simeq H^1(X,{\mathfrak t}_{\R}^*)^*/H_1(X,\Lambda_*(T)),
\end{gather}
where ${\mathfrak t}_{\R} = \R \underset{\Z}\times \Lambda_*(T)$ is the
split real form of the complex Lie algebra~${\mathfrak t}$ of~$T$.

As in Section \ref{higher}, for any
\begin{gather*}
\ga \in H^1(X,\Lambda^*(T)),
\end{gather*}
the image of $\ga$ in $H^1(X,{\mathfrak t}_{\R}^*)$ is represented by
a unique ${\mathfrak t}_{\R}^*$-valued one-form on $X$ that may be
written as
\begin{gather*}
\omega_\ga + \ol\omega_\ga,
\end{gather*}
where $\omega_\ga \in H^0\big(X,\Omega^{1,0}\big) \underset{\C}\otimes
{\mathfrak t}^*$ is a holomorphic ${\mathfrak t}^*$-valued one-form.

On the other hand, the image of $\ga$ in $H^1(X,{\mathfrak t}_{\R}^*)$
gives rise to a linear functional
\begin{gather*}
\varphi_\ga\colon \ H^1(X,{\mathfrak t}_{\R}^*)^* \to \R
\end{gather*}
satisfying $\varphi_\ga(\beta) \in \Z$ for all $\beta \in
H_1(X,\Lambda_*(T))$. Therefore, according to formula \eqref{JacT},
${\rm e}^{2\pi{\rm i} \varphi_{\ga}}$ is a well-defined function on
$\Bun^0_T(X)$. These are the Fourier harmonics on $\Bun^0_T(X)$.

In the same way as in Section~\ref{higher}, we prove the following result.

\begin{thm} \label{XT eig}
The functions ${\rm e}^{2\pi{\rm i} \varphi_{\ga}}$, $\ga \in
H^1(X,\Lambda^*(T))$, form an orthogonal basis of joint eigenfunctions
of the Hecke operators $_{p_0}H^{\cmu}_p$, $p \in X$, $\cmu \in
\Lambda_*(T)$, on $L^2\big(\Bun^0_T(X)\big)$. The eigenvalues of
$_{p_0}H^{\cmu}_p$ are given by the right hand side of the formula
\begin{gather} \label{eigomega1}
_{p_0}H^{\cmu}_p \cdot {\rm e}^{2\pi{\rm i} \varphi_{\ga}} = \cmu \left( {\rm e}^{2 \pi
 {\rm i} \int_{p_0}^p (\omega_\ga + \ol{\omega}_\ga)} \right) \; {\rm e}^{2\pi
 {\rm i} \varphi_{\ga}}.
\end{gather}
\end{thm}
Let us explain the notation we used on the right hand side of formula~\eqref{eigomega1}: denote by $\LT$ the Langlands dual torus to~$T$. We have $\Lambda_*\big(\LT\big) = \Lambda^*(T)$ and $\Lambda^*\big(\LT\big) =
\Lambda_*(T)$. The eigenvalue of the Hecke operator
$_{p_0}H^{\cmu}_p$, $p \in X$, $\cmu \in \Lambda^*\big(\LT\big)$, on the function
${\rm e}^{2\pi{\rm i} \varphi_{\ga}}$ is equal to the value of the
character $\cmu$ of~$\LT$ on the $\LT$-valued function $F_\gamma$ on~$X$
\begin{gather*} 
F_\ga(p) = {\rm e}^{2 \pi{\rm i} \int_{p_0}^p (\omega_\ga + \ol{\omega}_\ga)}, \qquad \ga
\in H^1(X,\Lambda^*(T)) = H^1\big(X,\Lambda_*\big(\LT\big)\big).
\end{gather*}
This function actually takes values in the compact form~$\LT_u$ of
$\LT$ and may be interpreted as the holonomy of the unitary connection
\begin{gather} \label{nablaga}
\nabla_\gamma = d - 2\pi{\rm i} (\omega_\ga + \ol{\omega}_\ga)
\end{gather}
on the trivial $\LT_u$-bundle on $X$ over (no matter which) path from
$p_0$ to $p$. As in the case of $T={\mathbb G}_m$, each of these
connections has trivial monodromy.

\section{Non-abelian case} \label{nonabelian}

In this section we try to generalize to the case of a non-abelian
group $G$ the results obtained in the previous section for abelian~$G$.

\subsection[Spherical Hecke algebra for groups over $\Fq\ppart$]{Spherical Hecke algebra for groups over $\boldsymbol{\Fq\ppart}$}

In the case of the function field of a curve $X$ over a finite field,
the Hecke operators attached to a closed point $x$ of $X$ generate the
spherical Hecke algebra ${\mc H}(G(\Fq\ppart),G(\Fq[[t]]))$. As a~vector space, it is the space of $\C$-valued functions on the group
$G(\Fq\ppart)$ that are bi-invariant with respect to the subgroup
$G(\Fq[[t]])$ (here $\Fq$ is the residue field of~$x$). This vector
space is endowed with the convolution product defined by the formula
\begin{gather*} 
(f_1 \star f_2)(g) = \int f_1\big(gh^{-1}\big) f_2(h)\,{\rm d}h,
\end{gather*}
where ${\rm d}h$ stands for the Haar measure on $G(\Fq\ppart)$ normalized so
that the volume of the subgroup $G(\Fq[[t]])$ is equal to~$1$ (in this
normalization, the characteristic function of~$G(\Fq[[t]])$ is the
unit element of the convolution algebra). The Haar measure can be
defined because~$G(\Fq\ppart)$ is a~{\em locally compact} group.

The resulting convolution algebra ${\mc H}(G(\Fq\ppart),G(\Fq[[t]]))$
is commutative and we have the Satake isomorphism between this algebra
and the complexified representation ring~$\on{Rep} \LG$ of the
Langlands dual group~$\LG$.

\subsection[Is there a spherical Hecke algebra for groups over $\C\ppart$?]{Is there a spherical Hecke algebra for groups over $\boldsymbol{\C\ppart}$?} \label{spher hecke}

In contrast to the group $G(\Fq\ppart)$, the group $G(\C\ppart)$ is
{\em not} locally compact. Therefore it does not carry a Haar
measure. Indeed, the field $\C\ppart$ is an example of a~two-dimensional local field, in the terminology of~\cite{Fesenko1},
more akin to~$\Fq\zpart\ppart$ or~${\mathbb Q}_p\ppart$ than to
$\Fq\ppart$ or~${\mathbb Q}_p$.

Ivan Fesenko has developed integration theory for the two-dimensional
local fields \cite{Fesenko1,Fesenko2}, and his students have extended
it to algebraic groups over such fields \cite{Morrow2,Morrow1,Waller},
but this theory is quite different from the familiar case of~$G(\Fq\ppart)$.

First, integrals over $\C\ppart$ and $G(\C\ppart)$ take values not in
real numbers, but in formal Laurent power series ${\mathbb
 R}(\!(X)\!)$, where $X$ is a formal variable. Under certain
restrictions, the value of the integral is a polynomial in $X$; if so,
then one could set $X$ to be equal to a real number. But this way one
might lose some important properties of the integration that we
normally take for granted.

Second, if $S$ is a Lebesque measurable subset of $\C$, then according
to \cite{Fesenko1,Fesenko2}, the measure of a subset of $\C\ppart$ of
the form
\begin{gather*} 
S t^i + t^{i+1}\C[[t]],
\end{gather*}
is equal to $\mu(S) X^i$, where $\mu(S)$ is the usual Lebesque measure
of $S$. In particular, this means that the measure of the subset~$\C[[t]]$ of $\C\ppart$ is equal to 0, as is the measure of the subset~$t^n\C[[t]]$ for any $n \in \Z$. Contrast this with the fact that
under a suitably normalized Haar measure on~$\Fq\ppart$, the measure
of $t^n \Fq[[t]]$ is equal to~$q^{-n}$. Thus, if we take as~$G$ the
additive group, it's not even clear how to define a unit element in
the would-be spherical Hecke algebra (which would be the
characteristic function of the subset $\Fq[[t]]$ in the case of the
field~$\Fq\ppart$). The situation is similar in the case of a general group~$G$.

For this reason, according to Waller~\cite{Waller:let}, from the point
of view of the two-dimensional integration theory it would make more
sense to consider distributions on $G\ppart$ that are bi-invariant not
with respect to $G(\C[[t]])$, but its subgroup~$\wt{K}$ consisting of
those elements $g(t) \in G(\C[[t]])$ for which~$g(0)$ belongs to a
compact subgroup~$K$ of~$G(\C)$. This would be similar to the
construction used in representation theory of complex Lie groups,
where one considers, for example, the space of distributions on the
group $G(\C)$ supported on a compact subgroup~$K$ with a~natural
convolution product~\cite{KV}.\footnote{I thank David Vogan for
 telling me about this construction, and the reference.} For
instance, if $K= \{1 \}$, the resulting algebra is~$U(\g)$, the
universal enveloping algebra of the Lie algebra~$\g$ of~$G(\C)$ (see
also Remark~\ref{nf} below.)

Perhaps, a convolution product on some space of distributions of this
kind can be defined for $G(\C\ppart)$, but from the structure of the
double cosets of $\wt{K}$ in $G(\C\ppart)$ it is clear that this
algebra would not be pertinent to defining Hecke operators on
$\Bun_G$.

Another option is to consider motivic integration theory. A motivic
version of the Haar measure has in fact been defined by Julia Gordon
\cite{Gordon} for the group $G(\C\ppart)$ (it may also be obtained in
the framework of the general theories of Cluckers--Loeser \cite{CL} or
Hrushovski--Kazhdan \cite{HK1}). Also, in a recent paper \cite{CCH} it
was shown that the spherical Hecke algebra ${\mc
 H}(G(\Fq\ppart),G(\Fq[[t]]))$ can be obtained by a certain
specialization from its version in which the ordinary integration with
respect to the Haar measure on $G(\Fq\ppart)$ is replaced by the
motivic integration with respect to the motivic Haar
measure. Presumably, one could carry some of the results of \cite{CCH}
over to the case of $\C\ppart$.

However, this does not seem to give us much help, for the following
reason: the motivic integrals over a ground field $k$ take values in a
certain algebra ${\mathcal M}_k$, which is roughly speaking a
localization of the Grothendieck ring of algebraic varieties over
$k$. In the case of the ground field $\Fq$, the algebra ${\mathcal
 M}_{\Fq}$ is rich, and ordinary integrals may be recovered from the
motivic ones by taking a homomorphism from ${\mathcal M}_{\Fq}$ to
$\R$ sending the class of the affine line over $\Fq$ to $q$. But in
the case of the ground field $\C$, the structure of the algebra ${\mc
 M}_{\C}$ appears to be very different (for example, it has divisors
of zero), and this construction does not work. In fact, it seems
that there are very few (if any) known homomorphisms from ${\mc
 M}_{\C}$ to positive real numbers, besides the Euler
characteristic.\footnote{I learned this from David Kazhdan (private
 communication).} Perhaps, taking the Euler characteristic, one can
obtain a non-trivial convolution algebra structure on the space of
$G(\C[[t]])$ bi-invariant functions on $G(\ppart)$ (one may wonder
whether it could be interpreted as a kind of $q \to 1$ limit of ${\mc
 H}(G(\Fq\ppart),G(\Fq[[t]]))$), but it is doubtful that this
algebra could be useful in any way for defining an analytic theory of
automorphic forms on $\Bun_G$ for complex algebraic curves.

\subsection{An attempt to define Hecke operators} \label{attempt}

What we are interested in here, however, is not the spherical Hecke
algebra itself, but rather the action of the corresponding Hecke
operators on automorphic functions. After all, in Section~\ref{abelian} we were able to define Hecke operators without any
reference to a convolution algebra on the group~$\C\ppart^\times$. The
abelian case, however, is an exception in that the action of the Hecke
operators did not require integration. In the non-abelian case (in
fact, already for $G={\rm GL}_2$), integration is necessary, and this
presents various difficulties, which I illustrate below with some
concrete examples in the case of ${\rm GL}_2$ and an elliptic curve.

Recall that for curves over $\Fq$, the unramified automorphic
functions are functions on the double quotient
\begin{gather} \label{double q}
G(F)\bs G({\mathbb A}_F)/G(\OO_F),
\end{gather}
where $F=\Fq(X)$, and $X$ is a curve over $\Fq$. The action of Hecke
operators on functions on this double quotient can be defined by means
of certain correspondences, and we can try to imitate this definition
for complex curves.

To this end, we take the same double quotient \eqref{double q} with
$F=\C(X)$, where $X$ is a curve over~$\C$. As in the case of $\Fq$,
this is the set of equivalence classes of principle $G$-bundles on
$X$. The Hecke correspondences can be conveniently defined in these
terms.

For instance, consider the case of ${\rm GL}_2$ and the first Hecke operator
(for a survey of the general case, see, e.g.,~\cite[Section~3.7]{F:rev}). Then we have the Hecke correspondence ${\mc H}{\rm ecke}_{1,x}$,
where~$x$ is a~closed point of~$X$:
\begin{gather} \label{Hecke cor}
\begin{array}{@{}ccccc}
& & {\mc H}{\rm ecke}_{1,x} & & \\
& \stackrel{\hl}\swarrow & & \stackrel{\hr}\searrow & \\
\Bun_{{\rm GL}_2} & & & & \Bun_{{\rm GL}_2}.
\end{array}
\end{gather}
Here ${\mc H}{\rm ecke}_{1,x}$ is the moduli stack classifying the
quadruples
\begin{gather*}(\M,\M',\beta\colon \M'\hookrightarrow\M),\end{gather*} where $\M$ and $\M'$ are
points of~$\Bun_{{\rm GL}_2}$, which means that they are rank two vector
bundles on~$X$, and~$\beta$ is an embedding of their sheaves of
(holomorphic) sections $\beta\colon \M'\hookrightarrow\M$ such that~$\M/\M'$
is supported at~$x$ and is isomorphic to the skyscraper sheaf $\OO_x =
\OO_X/\OO_X(-x)$. The maps are defined by the formulas
$\hl(\M,\M',\beta)=\M$, $\hr(\M,\M',\beta)=\M'$.

It follows that the points of the fiber of ${\mc H}{\rm ecke}_{1,x}$ over
${\mc M}$ in the ``left'' $\Bun_{{\rm GL}_2}$ correspond to all locally free
subsheaves ${\mc M}' \subset {\mc M}$ such that the quotient~${\mc
 M}/{\mc M}'$ is the skyscraper sheaf~$\OO_x$. Defining such~${\mc
 M}'$ is the same as choosing a line~$L$ in the dual space~${\mc
 M}_x^*$ to the fiber of~${\mc M}$ at~$x$ (which is a two-dimensional
complex vector space). The sections of the corresponding sheaf~${\mc
 M}'$ (over a Zariski open subset of~$X$) are the sections of ${\mc
 M}$ that vanish along~$L$, i.e., sections $s$ which satisfy the
equation $\langle v,s(x) \rangle = 0$ for a non-zero $v \in L$.

Thus, the fiber of ${\mc H}{\rm ecke}_{1,x}$ over ${\mc M}$ is isomorphic to
the projectivization of the two-dimen\-sio\-nal vector space ${\mc
 M}_x^*$, i.e., to $\C{\mathbb P}^1$. We conclude that ${\mc H}{\rm ecke}_{1,x}$ is a $\C{\mathbb P}^1$-fibration over the ``left'' $\Bun_{{\rm GL}_2}$ in the diagram~\eqref{Hecke cor}. Likewise, we
obtain that ${\mc H}{\rm ecke}_{1,x}$ is a~$\C{\mathbb P}^1$-fibration over
the ``right'' $\Bun_{{\rm GL}_2}$ in~\eqref{Hecke cor}.

In the geometric theory, we use the correspondence \eqref{Hecke cor}
to define a {\em Hecke functor} $\on{H}_{1,x}$ on the (derived)
category of $D$-modules on $\Bun_{{\rm GL}_2}$:
\begin{gather*} 
\on{H}_{1,x}(\K) = \hl{}_* \hr^*(\K)[1].
\end{gather*}
A $D$-module $\K$ is called a Hecke eigensheaf if we have
isomorphisms
\begin{gather} \label{eigen-property}
\imath_{1,x}\colon \ \on{H}_{1,x}(\K) \overset{\sim}\longrightarrow \C^2
\boxtimes \K \simeq \K \oplus \K, \qquad \forall\, x \in X,
\end{gather}
and in addition have similar isomorphisms for the second set of Hecke
functors $\on{H}_{2,x}$, $x \in X$. These are defined similarly to the
Hecke operators for ${\rm GL}_1$, as the pull-backs with respect to the
morphisms sending a rank two bundle $\M$ to $\M(x)$ (they change the
degree of~$\M$ by~$2$). (As explained in \cite[Section~1.1]{FGV}, the
second Hecke eigensheaf property follows from the first together with
a certain $S_2$-equivariance condition.)

Thus, if $\K$ is a Hecke eigensheaf, we obtain a family of
isomorphisms \eqref{eigen-property} for all $x \in X$, and similarly
for the second set of Hecke functors. We then impose a stronger
requirement that the two-dimensional vector spaces appearing on the
right hand side of~\eqref{eigen-property} as ``eigenvalues'' fit
together as stalks of a single rank two local system~${\mc E}$ on
$X$ (and similarly for the second set of Hecke functors, where the
eigenvalues should be the stalks of the rank one local system~$\wedge^2 {\mc E}$ on~$X$; this is, however, automatic if we impose the $S_2$-equivariance condition from \cite[Section~1.1]{FGV}). If
that's the case, we say that $\K$ is a {\em Hecke eigensheaf} with the
{\em eigenvalue} ${\mc E}$. This is explained in more detail, e.g., in
\cite[Section~3.8]{F:rev}.

The first task of the geometric theory (in the case of $G={\rm GL}_2$) is to
show that such a~Hecke eigensheaf on $\Bun_{{\rm GL}_2}$ exists for every
irreducible rank two local system~${\mc E}$ on~$X$. This was
accomplished by Drinfeld in~\cite{Dr1}, a groundbreaking work that
was the starting point of the geometric theory. We now know that the
same is true for $G={\rm GL}_n$~\cite{FGV,Ga1} and in many other cases.

Now let's try to adapt the diagram~\eqref{Hecke cor} to
functions. Thus, given a function~$f$ on the set of $\C$-points, we
wish to define the action of the first Hecke operator~$H_{1,x}$ on it
by the formula
\begin{gather} \label{formula He1}
(H_{1,x} \cdot f)(\M) = \underset{\M' \in \hl^{-1}(\M)}\int f(\M') \, {\rm d}\M'.
\end{gather}
Thus, we see that the result must be an integral over the complex
projective line~$\hl^{-1}(\M)$. The key question is: what is the
measure~${\rm d}\M'$?

Herein lies a crucial difference with the abelian case considered in
Section \ref{abelian}: in the abelian case every Hecke operator acted
by pull-back of a function, so no integration was needed. But in the
non-abelian case, already for the first Hecke operators $H_{1,x}$ in
the case of $G={\rm GL}_2$, we must integrate functions over the projective
lines $\hl^{-1}(\M)$, where $\M \in \Bun_{{\rm GL}_2}(\C)$.

Note that if our curve were over a finite field, this integration is
in fact a summation over a~finite set of $q+1$ elements, the number of
points of $\pone$ over~$\Fq$, where~$\Fq$ is the residue field of the
closed point~$x$ at which we take the Hecke operator. The terms of
this summation correspond to points of the fibers
$\hl^{-1}(\M)$. Being finite sums, these integrals are always
well-defined if our curve is over~$\Fq$. For curves over~$\C$, this is
not so, and this creates major problems, as we will see below.

\subsection{The case of an elliptic curve} \label{attempt ell}

Let's look at the case of an elliptic curve $X$. If it is defined over
a finite field, the fibers $\hl^{-1}(\M)$ appearing in the Hecke
operators have been described explicitly in
\cite{Alvarenga, Lorscheid}, using the classification of rank two
bundles on elliptic curves due to Atiyah~\cite{Atiyah}.

For a complex elliptic curve~$X$, the fibers $\hl^{-1}(\M)$ have been
described explicitly in~\cite{Boozer}. In~\cite{L:analyt}, Langlands
attempts to describe them in the language of adeles, which is more
unwieldy than the vector bundle language used in \cite{Boozer} and
hence more prone to errors. As the result of his computations,
Langlands states on p.~18 of~\cite{L:analyt}:

``The dimension $\dim(g\Delta_1/G(\OO_x))$ [which
is our $\hl^{-1}(\M)$ if $\M$ is the rank two bundle corresponding to
the adele $g \in {\rm GL}_2({\mathbb A}_F)$] is always equal to~0 \dots\ Hence
the domain of integration in [adelic version of our formula~\eqref{formula He1} above] is a finite set''.

This statement is incorrect. First of all, as we show below, there are
rank two bundles $\M$ on an elliptic curve for which there are {\em
 infinitely many} non-isomorphic bundles in the fiber~$\hl^{-1}(\M)$
(in fact, we have a continuous family of non-isomorphic bundles
parametrized by the points of~$\hl^{-1}(\M)$). Second, even if there
are finitely many isomorphism classes among those $\M'$ which appear
in the fiber $\hl^{-1}(\M)$ for a fixed~$\M$, this does not mean that
we are integrating over a~finite set.

In fact, according to formula~\eqref{formula He1} (whose adelic
version is formula~(10) of~\cite{L:analyt}), for every~$\M$, the fiber
$\hl^{-1}(\M)$ of the Hecke correspondence over which we are supposed
to integrate is {\em always isomorphic to} $\C\pone$ if we take into
account the automorphism groups of the bundles involved. In the adelic
language, the automorphism group $\on{Aut}({\mc M})$ of a bundle ${\mc
 M}$ may be described, up to an isomorphism, as follows: $\M$
corresponds to a point in the double quotient~\eqref{double q}; we
lift this point to $G(F) \bs G({\mathbb A}_F)$ and take its stabilizer
subgroup in~$G(\OO_F)$.

The necessity to take into account these automorphism groups is
well-know in the case of curves defined over $\Fq$. In this case, the
measure on the double quotient~\eqref{double q} induced by the
Tamagawa measure assigns (up to an overall factor) to a point not~$1$
but $1/|\on{Aut}({\mc M})|$ (this measure is well-defined if we work
over $\Fq$ because then the group $\on{Aut}({\mc M})$ is finite for
any ${\mc M}$; however, this is not so over $\C$). With respect to
this correctly defined measure, the fiber $\hl^{-1}(\M)$ for any ${\mc
 M}$ and any $\Fq$-point~$x$ of~$X$ can be identified with the set of
$\Fq$-points of the projective line over~$\Fq$, with each point having
measure~$1$.

As a concrete illustration, consider the following example.

\begin{Example} Let $\M = {\mc L}_1
\oplus {\mc L}_2$, where ${\mc L}_1$ and ${\mc L}_2$ are two line
bundles of degrees $d_1$ and $d_2$ such that $d_1>d_2+1$. Then the
vector bundles $\M'$ that appear in the fiber $\hl^{-1}(\M)$ are
isomorphic to either $\M'_1 = {\mc L}_1(-x) \oplus {\mc L}_2$ or
$\M'_2 = {\mc L}_1 \oplus {\mc L}_2(-x)$. However, the groups of
automorphisms of these bundles are different: each of them is a
semi-direct product of the group $(\C^\times)^2$ of rescalings of the
two line bundles appearing in a direct sum decomposition and an
additive group, which is $\on{Hom}({\mc L}_2,{\mc L}_1)$ for $\M$;
$\on{Hom}({\mc L}_2,{\mc L}_1(-x))$ for $\M'_1$; and $\on{Hom}({\mc
 L}_2(-x),{\mc L}_1)$ for $\M'_2$.

Under our assumption that $d_1>d_2+1$, we find that the latter groups
are isomorphic to $\C^{d_1-d_2-1}$, $\C^{d_1-d_2-2}$, and
$\C^{d_1-d_2}$, respectively. Thus, the automorphism group of $\M'_1$
is ``smaller'' by one copy of the additive group ${\mathbb G}_a$ than
that of $\M$, whereas the automorphism group of $\M'_2$ is larger than
that of $\M$ by the same amount.

This implies that the fiber $\hl^{-1}(\M)$ is the union of a complex
affine line worth of points corresponding to $\M'_1$ and a single
point corresponding to $\M'_2$.

If we worked over $\Fq$, we would find that
\begin{gather*}
(H_{1,x} \cdot f)(\M) = q f(\M'_1) + f(\M'_2),
\end{gather*}
with the factor of $q$ representing the number of points of the affine
line (see \cite{Alvarenga, Lorscheid}). Over~$\C$, we formally obtain
the sum of two terms: (1)~an {\em integral} of the constant function
taking value~$f(\M'_1)$ on an open dense subset of~$\C\pone$
isomorphic to the affine line, and (2)~a single term~$f(\M'_2)$
corresponding to the remaining point. Is there an integration measure
of $\C\pone$ that would render this sum meaningful?

If we use a standard integration measure on $\C\pone$, then the answer
would be $f(\M'_1)$ multiplied by the measure of the affine line. The
second term would drop out, as it would correspond to a subset
(namely, a point) of measure zero. If we want to include the second
bundle (which we certainly do for Hecke operators to be meaningful),
then the measure of this point has to be non-zero. But we also expect
our measure on~$\C\pone$ to be invariant (indeed, we cannot {\em a~priori} distinguish a special point on each of these projective
lines). Therefore the measure of every point of~$\C\pone$ would have
to be given by the same non-zero number. But then our integral would
diverge. It is not clear how one could regularize these divergent
integrals in a uniform and meaningful way.
\end{Example}

Next, we give an example in which the fiber $\hl^{-1}(\M)$ is a~continuous family of {\em non-isomorphic} vector bundles (thus directly contradicting the above statement from~\cite{L:analyt}).

\begin{Example}\label{Example2} Let $\M$ be the indecomposable rank two degree~1 vector bundle~$F_2(x)$ (in the notation of~\cite{L:analyt}) which is a unique, up to an isomorphism, non-trivial extension
\begin{gather*} 
0 \to \OO_X \to F_2(x) \to \OO_X(x) \to 0.
\end{gather*}
In this case, as shown in \cite[Section~4.3]{Boozer}, the fiber
$\hl^{-1}(F_2(x))$ may be described in terms of a~canonical
two-sheeted covering $\pi\colon \Pic^0(X) \to \C\pone = \hl^{-1}(F_2(x))$
ramified at 4~points such that (1)~if $a \in \hl^{-1}(F_2(x))$ is
outside of the ramification locus, then $\pi^{-1}(a) = \big\{ {\mc L}[a],
{\mc L}[a]^{-1} \big\}$, where ${\mc L}[a]$ is a degree~$0$ line bundle on~$X$; and (2)~the fibers over the 4~ramification points are the four square roots ${\mc L}_i$, $i=1,\ldots,4$, of the trivial line bundle on~$X$.

Namely, the vector bundle $\M'(a)$ corresponding to a point $a \in
\hl^{-1}(F_2(x))$ is described in terms of~$\pi$ as follows (note that
in~\cite{Boozer} the bundle $F_2(x)$ is denoted by~$G_2(x)$ and ${\mc
 L}[a]$ is denoted by~${\mc L}(a)$):
\begin{itemize}\itemsep=0pt
\item if $a \in \hl^{-1}(F_2(x))$ is outside of the ramification locus,
 then $\M'(a) = {\mc L}[a] \oplus {\mc L}[a]^{-1}$;

\item if $a$ is a ramification point corresponding to the line bundle
 ${\mc L}_i$, then $\M'(a) = {\mc L}_i \otimes F_2$, where $F_2$ is
 the unique, up to an isomorphism, non-trivial extension of $\OO_X$
 by itself.
\end{itemize}

According to the Atiyah's classification, the bundles $\M'(a)$ and
$\M'(b)$ corresponding to different points $a \neq b$ in
$\hl^{-1}(F_2(x))$ are non-isomorphic. Thus, there is an infinite
continuous family of non-isomorphic vector bundles appearing in the
fiber $\hl^{-1}(F_2(x))$ in this case.

One gets a similar answer for $\M=F_2(x) \otimes {\mc L}$, where ${\mc
 L}$ is an arbitrary line bundle on $X$ (note that unlike the vector
bundles discussed in the previous example, all of the bundles $F_2(x)
\otimes {\mc L}$ are {\em stable}). This means that the value of
$H_{1,x} \cdot f$ at bundles $\M$ of this form depends on the choice
of a measure of integration on $\hl^{-1}(\M)$.
\end{Example}

It is not clear whether it is possible to define these measures for
different $\M$ and different $x \in X$ in a consistent and meaningful
way, so that they would not only yield well-defined integrals but that
the corresponding operators $H_{1,x}$, $x \in X$, would commute with
each other and with the second set of Hecke operators $H_{2,x}$, $x \in X$.

In \cite{L:analyt}, Langlands sidesteps these problems and instead
defines his versions of the Hecke operators by explicit formulas. Let
$D$ (resp. $U$) be the substacks of $\Bun_{{\rm GL}_2}$ of an elliptic curve
parametrizing rank two vector bundles on an elliptic curve $X$ that
are decomposable (resp. indecomposable) as direct sums of line
bundles. Points in both substacks can be explicitly described using
Atiyah's classification results~\cite{Atiyah}. If we forget the
automorphism groups of these rank two bundles, we obtain algebraic
varieties ${\mathfrak D}$ and ${\mathfrak U}$, the former isomorphic
to $\on{Sym}^2(\Pic(X))$ and the latter isomorphic to a disjoint union
of $\Pic(X)$ and $\Pic(X)/\Pic_2(X)$, where $\Pic_2(X)$ is the
subgroup of line bundles ${\mc L}$ such that ${\mc L}^{\otimes 2}
\simeq \OO_X$ (they correspond to the indecomposable bundles of even
and odd degrees, respectively).

In \cite[pp.~20--21]{L:analyt}, Langlands defines his versions of the
Hecke operators as linear operators acting on the direct sum
$L^2({\mathfrak D}) \oplus L^2({\mathfrak U})$. He postulates (as he
writes on p.~21, ``by decree!'') 
that these Hecke operators should act on
this space diagonally, i.e., preserving each of the two subspaces
$L^2({\mathfrak D})$ and $L^2({\mathfrak U})$.

However, the idea of treating the moduli stack of rank two bundles on
$X$ as the disjoint union of the varieties ${\mathfrak D}$ and
${\mathfrak U}$ sounds problematic. Indeed, in the moduli stack
$\Bun_{{\rm GL}_2}$ the substacks $D$ and $U$ are ``glued'' together in a
non-trivial way.\footnote{Note that if we were to consider instead the
 moduli space of semi-stable bundles, then, depending on the
 stability condition we choose, some of the decomposable bundles in
 ${\mathfrak D}$ would have to be removed, or identified with the
 indecomposable ones in ${\mathfrak U}$. Considering the moduli space
 of semi-stable bundles is, however, problematic for a different
 reason: it is not preserved by the Hecke correspondences.} If we
tear them apart, we are at the same time tearing apart the projective
lines $\hl^{-1}(\M)$ appearing as the fibers of the Hecke
correspondences. The idea of defining the Hecke operators in such away
that we keep the part corresponding to the bundles of one type and
throw away the part corresponding to the bundles of the other type
sounds even more problematic, for the following reason.

Consider the stable indecomposable bundle $F_2(x)$ described in
Example~\ref{Example2} above. It corresponds to a point of ${\mathfrak U}$. But, as
explained in Example~\ref{Example2}, all but 4 points in the Hecke fiber
$\hl^{-1}(F_2(x)) \simeq \C\pone$ of $F_2(x)$ correspond to
semi-stable {\em decomposable} bundles, i.e., points of ${\mathfrak
 D}$. (Furthermore, the bundles corresponding to different points in
$\hl^{-1}(F_2(x))$ are not isomorphic to each other.) It is not clear
why one would want to throw away this {\em open dense subset} of the
Hecke fiber. More examples of this nature are given below.

\begin{Example} \label{Example3} Let $\M = \OO_X \oplus \OO_X(x)$. Then, as
explained in~\cite{Boozer}, there are two points in the fiber
$\hl^{-1}(\M)$, corresponding to $\M'_1 = \OO_X \oplus \OO_X$ and
$\M'_2 = \OO_X(-x) \oplus \OO_X(x)$, and each point in the complement
(which is isomorphic to $\C^\times$) corresponds to the indecomposable
bundle $F_2$. Thus, we see that an open dense subset of the fiber
$\hl^{-1}(\M)$ of the Hecke correspondence over a decomposable rank
two bundle $\M$ (corresponding to a point in ${\mathfrak D}$) consists
of {\em indecomposable} bundles (corresponding to points of ${\mathfrak
 U}$).

Now, if we were to treat $\M'_1 = \OO_X \oplus \OO_X$ and $\M'_2 =
\OO_X(-x) \oplus \OO_X(x)$ as belonging to a~different connected
component of $\Bun_{{\rm GL}_2}$ than~$F_2$, then what to make of the
integral~\eqref{formula He1}? It would seemingly break into the sum of
two points and an integral over their complement. That would be fine
in the case of a curve over~$\Fq$: we would simply obtain the formula
\begin{gather*}
(H_{1,x} \cdot f)(\M) = f(\M'_1) + f(\M'_2) + (q-1) f(F_2),
\end{gather*}
with the factor of $(q-1)$ being the number of points of $\pone$
without two points (see \cite{Alvarenga, Lorscheid}). But over complex
numbers we have to integrate over $\C^\times$. We would therefore have
to somehow combine summation over two points and integration over
their complement. As in another example of this nature that we
considered above, it is not clear that there exists an integration
measure that would achieve this in a consistent and meaningful
fashion.
\end{Example}

\begin{Example}\label{Example4} Similarly, if $\M = {\mc L}_1 \oplus {\mc
 L}_2$, where ${\mc L}_1$ and ${\mc L}_2$ are non-isomorphic line
bundles of degree $0$, then the fiber of the Hecke correspondence
$\hl^{-1}(\M)$ over $\M$ has two points corresponding $\M'_1 = {\mc
 L}_1(-x) \oplus {\mc L}_2$ and $\M'_2 = {\mc L}_1 \oplus {\mc
 L}_2(-x)$, and every point in the complement of these two points
(which is isomorphic to $\C^\times$) corresponds to the indecomposable
vector bundle $F_2(x) \otimes {\mc L}(-x)$, where ${\mc L}^{\otimes 2}
\simeq {\mc L}_1 \otimes {\mc L}_2$. Again, we see that an open dense
subset of the fiber $\hl^{-1}(\M)$ of the Hecke correspondence over a
decomposable rank two bundle $\M$ (corresponding to a point in
${\mathfrak D}$) consists of {\em indecomposable bundles}
(corresponding to points of ${\mathfrak U}$). As in the previous
example, it is not clear how to integrate over $\hl^{-1}(\M)$.
\end{Example}

Incidentally, Example~\ref{Example4} shows that for any $x \in X$ and any pair
of degree zero line bundles~${\mc L}$,~${\mc L}_1$ on~$X$, there is a
continuous family of rank two vector bundles on $X$ over an affine
line ${\mathbb A}^1$ that are isomorphic to $F_2(x) \otimes {\mc L}$
away from $0 \in {\mathbb A}^1$ and to ${\mc L}_1(x) \oplus \big({\mc
 L}^{\otimes 2} \otimes {\mc L}_1^{-1}\big)$ at the point $0 \in {\mathbb
 A}^1$. Likewise, Example~\ref{Example3} shows that there exist
continuous families of rank two vector bundles on $X$ over an affine
line ${\mathbb A}^1$ that are isomorphic to~$F_2$ away from $0 \in
{\mathbb A}^1$ and to $\OO_X \oplus \OO_X$ or to $\OO_X(-x) \oplus
\OO_X(x)$ at the point $0 \in {\mathbb A}^1$.

These examples illustrate the intricate (non-Hausdorff) topology of
$\Bun_{{\rm GL}_2}$ of an elliptic curve; in particular, the fact that the
substacks~$D$ and~$U$ of decomposable and indecomposable bundles are
glued together in a highly non-trivial fashion.

To summarize: in~\cite{L:analyt} Langlands defines his versions of
Hecke operators in the case of ${\rm GL}_2$ and an elliptic curve~$X$ in an
{\em ad hoc} fashion, without discussing the pertinent measures of
integration from the first principles. As far as I understand, his
definition is based on two assumptions: (i)~the
statement that the fibers of the Hecke correspondence are finite, see
the quote from~\cite{L:analyt} at the beginning of this subsection;
(ii)~postulating that the Hecke operators
should act on $L^2({\mathfrak D}) \oplus L^2({\mathfrak U})$
preserving each of the two direct summands. In this section, I have
outlined the issues with these
assumptions.\footnote{{\em Added in September 2019:} If $\Bun_G$
 contains an open dense substack of stable bundles, it is possible to
 define analogues of Hecke operators acting on compactly supported
 sections of the line bundle of half-densities on this substack (rather
 than functions), following the construction of A.~Braverman and
 D.~Kazhdan~\cite{BK1} in the non-archimedian case. The details will
 appear in~\cite{EFK1}. However, this construction cannot be applied in
 the case of elliptic curves because there is no such open dense
 substack in $\Bun_G$ in this case (unless we add some extra structures
 to $G$-bundles, such as parabolic structures).

 In a special case ($G={\rm PGL}_2, X=\pone$ with parabolic structures at
 four points~-- in this case $\Bun_G$ contains an open dense substack
 of stable bundles), explicit formulas for the Hecke operators were
 proposed by M.~Kontsevich \cite[Section~2.4]{Kon}.}

\section{An alternative proposal} \label{alt}

There is however another approach to the analytic theory of
automorphic functions for complex curves, proposed in a joint work
with Pavel Etingof and David Kazhdan~\cite{EFK}. In this section I~outline this approach.

\subsection{A toy model} \label{toy}

It is instructive to consider first a toy model for the questions we have been discussing. Over~$\Fq$, there is a well-understood finite-dimensional analogue of the spherical Hecke algebra of~$G(\Fq\ppart)$; namely, the Hecke algebra~${\mc H}_q(G)$ of~$B(\Fq)$ bi-invariant $\C$-valued functions on the group~$G(\Fq)$, where~$B$ is a Borel subgroup of a simple algebraic group~$G$.

As a vector space, this algebra has a basis labeled by the
characteristic functions~$c_w$ of the Bruhat--Schubert cells $B(\Fq) w
B(\Fq)$, where $w$ runs over the Weyl group of~$G$. The convolution
product on ${\mc H}_q(G)$ is defined using the constant measure
$\mu_q$ on the finite group $G(\Fq)$ normalized so that the measure of
$B(\Fq)$ is equal to $1$. Then the function $c_1$ is a unit element of
${\mc H}_q(G)$.

It is convenient to describe the convolution product on ${\mc H}_q(G)$
as follows: identify the $B$ bi-invariant functions on $G$ with
$B$-invariant functions on $G/B$ and then with $G$-invariant functions
on $(G/B) \times (G/B)$ (with respect to the diagonal action). Given
two $G$-invariant functions $f_1$ and $f_2$ on $(G/B) \times (G/B)$,
we define their convolution product by the formula
\begin{gather} \label{conv GB}
(f_1 \star f_2)(x,y) = \underset{G/B}\int f_1(x,z) f_2(z,y)\,{\rm d}z.
\end{gather}
Under this convolution product, the algebra ${\mc H}_q(G)$ is
generated by the functions~$c_{s_i}$, where the~$s_i$ are the simple
reflections in~$W$. They satisfy the well-known relations.

Observe also that the algebra ${\mc H}_q(G)$ naturally acts on the
space $\C[G(\Fq)/B(\Fq)]$ of $\C$-valued functions on~$G(\Fq)/B(\Fq)$. It acts on the right and commutes with the natural
left action of~$G(\Fq)$. Unlike the spherical Hecke algebra, ${\mc
 H}_q(G)$ is non-commutative. Nevertheless, we can use the
decomposition of the space $\C[G(\Fq)/B(\Fq)]$ into irreducible
representations of ${\mc H}_q(G)$ to describe it as direct sum of
irreducible representations of~$G(\Fq)$.

Now suppose that we wish to generalize this construction to the
complex case. Thus, we consider the group~$G(\C)$, its Borel subgroup
$B(\C)$, and the quotient $G(\C)/B(\C)$, which is the set of
$\C$-points of the flag variety $G/B$ over~$\C$. A naive analogue of
${\mc H}_q(G)$ would be the space~${\mc H}_{\C}(G)$ of~$B(\C)$
bi-invariant functions on~$G(\C)$. Therefore we have the following
analogues of the questions that we discussed above in the case of the
spherical Hecke algebra: Is it possible to define a measure of
integration on~$G(\C)$ that gives rise to a meaningful convolution
product on~${\mc H}_{\C}(G)$? Is it possible to use the resulting
algebra to decompose the space of $L^2$ functions on~$G(\C)/B(\C)$?

\looseness=-1 For example, consider the case of $G={\rm SL}_2$. Then $G/B = \pone$. The
Hecke algebra ${\mc H}_q({\rm SL}_2)$ has a basis consisting of two
elements, $c_1$ and $c_s$, which (in its realization as $G$-invariant
functions on $(G/B) \times (G/B)$ explained above) correspond to the
characteristic functions of the two ${\rm SL}_2$-orbits in $\pone \times
\pone$: the diagonal and its complement, respectively. Applying
formula \eqref{conv GB}, we obtain that
\begin{gather} \label{c1}
c_1 \star c_1 = c_1, \qquad c_1 \star c_s = c_s, \\ \label{cs}
c_s \star c_s = q c_1 + (q-1) c_s.
\end{gather}
The two formulas in \eqref{c1} follow from the fact that for each $x$
and $y$, in formula \eqref{conv GB} there is either a unique value of
$z$ for which the integrand is non-zero, or no such values. The
coefficients in formula \eqref{cs} have the following meaning: $q =
\mu_q\big({\mathbb A}^1\big)$, $q-1 = \mu_q\big({\mathbb A}^1\bs 0\big)$.

Now, if we try to adopt this to the case of $\pone$ over~$\C$, we
quickly run into trouble. Indeed, if we want~$c_1$ to be the unit
element, we want to keep the two formulas in~\eqref{c1}. But in order
to reproduce the second formula in~\eqref{c1}, we need a measure $dz$
on $\C\pone$ that would give us $\int \chi_u \, {\rm d}z = 1$ for every point
$u \in \C\pone$, where~$\chi_u$ is the characteristic function of~$u$. However, then the integral of this measure over the affine line
inside~$\C\pone$ would diverge, rendering the convolution product~$c_s
\star c_s$ meaningless.

Likewise, we run into trouble if we attempt to define an action of
${\mc H}_{\C}(G)$ on the space of functions on~$G(\C)/B(\C)$. Thus, we
see that the questions we asked above do not have satisfactory
answers, and the reasons for that are similar to those we discussed in
the previous section, concerning the spherical Hecke algebra and the
possibility of defining an action of Hecke operators on functions on~$\Bun_G$.

However, there are two natural variations of these questions that do
have satisfactory answers. The first possibility is to consider a {\em
 categorical} version of the Hecke algebra, i.e., instead of the
space of $B$-invariant functions, the category $D(G/B)^B$-mod of
$B$-equivariant $D$-modules on~$G/B$. According to a theorem of
Beilinson and Bernstein~\cite{BB}, we have an exact functor of global
sections (as $\OO$-modules) from this category to the category of
modules over the Lie algebra~$\g$ of~$G$, which is an equivalence with
the category of those $\g$-module which have a fixed character of the
center of~$U(\g)$ (the character of the trivial representation of
$\g$). This is the category that appears in the Kazhdan--Lusztig
theory, which gives rise, among other things, to character formulas
for irreducible $\g$-modules from the category ${\mc O}$. Furthermore,
instead of the convolution product on functions, we now have
convolution functors on a derived version of $D(G/B)^B$-mod. This is
the categorical Hecke algebra which has many applications. For
example, Beilinson and Bernstein have defined a categorical action of
this category on the derived category of the category $\OO$ (which may
be viewed as the category of $(\g,B)$ Harish-Chandra modules). This is
a special case of a rich theory.

Note that a closely related category of perverse sheaves may also
be defined over $\Fq$. Taking the traces of the Frobenius on the
stalks of those sheaves, we obtain the elements of the original Hecke
algebra~${\mc H}_q(G)$. This operation transforms convolution product
of sheaves into convolution product of functions. Thus, we see many
parallels with the geometric Langlands Program (for more on this, see
\cite[Section~1.3.3]{F:loop}). In particular, the spherical Hecke
algebra has a categorical analogue, for which a categorical version of
the Satake isomorphism has been proved~\cite{Ginzburg,Lusztig,MV}. In
other words, the path of categorification of the Hecke algebra ${\mc
 H}_q(G)$ is parallel to the path taken in the geometric Langlands
theory.

But there is also a second option: We can define a Hilbert space
$L^2(G(\C)/B(\C))$ as the completion of the space of half-densities on
$G(\C)/B(\C)$ with respect to the natural Hermitian inner product
(this has a generalization corresponding to twisting by line bundles
on $G/B$ as well as certain ``imaginary powers'' thereof). However,
instead of defining an action of a Hecke algebra on this Hilbert
space, one then uses a substitute: differential operators on $G/B$.

The Lie algebra $\g$ acts on $L^2(G(\C)/B(\C))$ by holomorphic vector
fields, and we have a commuting action of another copy of $\g$ by
anti-holomorphic vector fields. Therefore, the tensor product of two
copies of the center of $U(\g)$ acts by mutually commuting
differential operators.\footnote{The referee drew my attention to the
 fact that this is essentially the classical Gelfand--Naimark
 construction of the principal series representations of $G(\C)$ and
 the corresponding Harish-Chandra bimodules.} As we mentioned above,
both holomorphic and anti-holomorphic ones act according to the
central character of the trivial representation. However, the center
of $U(\g_c$), where $\g_c$ is a~compact form of the Lie algebra~$\g$,
also acts on~$L^2(G(\C)/B(\C))$ by commuting differential operators,
and this action is non-trivial. It includes the Laplace operator,
which corresponds to the Casimir element of~$U(\g_c)$.

We then ask what are the eigenfunctions and eigenvalues of these
commuting differential ope\-rators. This question has a meaningful
answer. Indeed, using the isomorphism $G/B \simeq G_c/T_c$, where~$T_c$ is a maximal torus of the compact form of~$G$, and the Peter--Weyl theorem, we obtain that~$L^2(G(\C)/B(\C))$ can be
decomposed as a direct sum of irreducible finite-dimensional
representations of~$\g_c$ which can be exponentiated to the group~$G_c$ of adjoint type, each irreducible representation~$V$ appearing
with multiplicity equal to the dimension of the weight~0 subspace~$V(0)$ of~$V$. Therefore the combined action of the center of
$U(\g_c)$ and the Cartan subalgebra ${\mathfrak t}_c$ of $T_c$ (acting
by vector fields) has as eigenspaces various weight components of
various irreducible representations $V$ of $\g_c$ tensored with~$V(0)$. All of these eigenspaces are finite-dimensional.

For instance, for $G={\rm SL}_2$ every eigenspace is one-dimensional, and so
we find that these differential operators have simple spectrum. In
fact, suitably normalized joint eigenfunctions of the center of
$U(\g_c)$ and~${\mathfrak t}_c$ are in this case the standard
spherical harmonics (note that in this case $G(\C)/B(\C) \simeq S^2$).

This discussion suggests we may be able to build a meaningful analytic
theory of automorphic forms on $\Bun_G$ if, rather than looking for
the eigenfunctions of Hecke operators (whose existence is questionable
in the non-abelian case, as we have seen), we look for the
eigenfunctions of a~commutative algebra of global differential
operators on~$\Bun_G$. It turns out that we are in luck: there exists
a large commutative algebra of differential operators acting on the
line bundle of half-densities on $\Bun_G$.

\begin{rem} \label{nf}The above discussion dovetails nicely with the intuition that comes
 from the theory of automorphic functions for a reductive group~$G$
 over a number field~$F$. Such a field has non-archimedian as well as
 archimedian completions. The representation theories of the
 corresponding groups, such as~$G(\Q_p)$ and~$G(\C)$, are known to follow different paths: for the
 former we have, in the unramified case, the spherical Hecke algebra
 and the Satake isomorphism. For the latter, instead of a~spherical
 Hecke algebra one usually considers the center of $U(\g)$ (or, more
 generally, the convolution algebra of distributions on $G(\C)$
 supported on its compact subgroup $K$, see \cite{KV}).

 Now let's replace a number field $F$ by a field of the form $F(X)$,
 where $X$ is a curve over~$F$. Then instead of the local fields
 $\Q_p$ we would have fields such as $\Q_p\ppart$, and instead of~$\C$ we would have~$\C\ppart$. In the former case we would have to
 consider the group~$G(\Q_p\ppart)$ and in the latter case, the group~$G(\C\ppart)$. For $G(\Q_p\ppart)$ there are meaningful analogues of
 the spherical Hecke algebra and the corresponding Satake
 isomorphism. They have been studied, in particular, in
 \cite{BK1,BK2,HK2, Kap,KL}. But in the case of $G(\C\ppart)$, just as
 in the case of a number field $F$ discussed above, it seems
 more prudent to consider the center of $U(\g\ppart)$ instead. As
 we show in the rest of this section, this approach leads to a rich
 and meaningful theory. Indeed, if we take the so-called critical
 central extension of~$\g\ppart$, then the corresponding completed
 enveloping algebra does contain a large center, as shown in
 \cite{FF} (see also \cite{F:wak,F:loop}). This center gives rise to
 a large algebra of global commuting differential operators on
 $\Bun_G$.
\end{rem}

\subsection[Global differential operators on $\Bun_G$]{Global differential operators on $\boldsymbol{\Bun_G}$} \label{global}

Let us assume for simplicity that $G$ is a connected,
simply-connected, simple algebraic group over $\C$. In~\cite{BD},
Beilinson and Drinfeld have described the algebra $D_G$ of global
holomorphic differential operators on $\Bun_G$ acting on the square
root $K^{1/2}$ of a canonical line bundle (which exists for any
reductive $G$ and is unique under our assumptions). They have proved
that $D_G$ is commutative and is isomorphic to the algebra of
functions on the space $\on{Op}_{\LG}(X)$ of $\LG$-opers on $X$. For a
survey of this construction and the definition of $\on{Op}_{\LG}(X)$,
see, e.g., \cite[Sections~8 and~9]{F:rev}. Under the above assumptions
on~$G$, the space~$\on{Op}_{\LG}(X)$ may be identified with the space
of all holomorphic connections on a particular holomorphic
$\LG$-bundle ${\mc F}_0$ on~$X$. In particular, it is an affine space
of dimension equal to $\dim \Bun_{\LG}$.

The construction of these global differential operators is similar to
the construction outlined in Section \ref{toy} above. Namely, they are
obtained in \cite{BD} from the central elements of the completed
enveloping algebra of the affine Kac--Moody algebra $\ghat$ at the
critical level, using the realization of $\Bun_G$ as a double quotient
of the formal loop group $G(\C\ppart)$ and the Beilinson--Bernstein
type localization functor. The critical level of $\ghat$ corresponds
to the square root of the canonical line bundle on $\Bun_G$. A theorem
of Feigin and myself \cite{FF} (see also \cite{F:wak,F:loop})
identifies the center of this enveloping algebra with the algebra of
functions on the space of $\LG$-opers on the formal punctured
disc. This is a local statement that Beilinson and Drinfeld use in the
proof of their theorem.

Now, we can use the same method to construct the algebra $\ol{D}_G$ of
global anti-holomorphic differential operators on $\Bun_G$ acting on
the square root $\ol{K}^{1/2}$ of the anti-canonical line bundle. The
theorem of Beilinson and Drinfeld implies that $\ol{D}_G$ is
isomorphic to the algebra of functions on the complex conjugate space
to the space of opers, which we denote by
$\ol{\on{Op}}_{\LG}(X)$. Under the above assumptions on $G$, it can be
identified with the space of all anti-holomorphic connections on the
$G$-bundle $\ol{\mc F}_0$ that is the complex conjugate of the
$G$-bundle ${\mc F}_0$. While ${\mc F}_0$ carries a~holomorphic
structure (i.e., a $(0,1)$-connection), $\ol{\mc F}_0$ carries a
$(1,0)$-connection (which one could call an ``anti-holomorphic
structure'' on~$\ol{\mc F}_0$). Just as a $(1,0)$, i.e., holomorphic,
connection on~${\mc F}_0$ completes its holomorphic structure to a
flat connection, so does a~$(0,1)$, i.e., anti-holomorphic, connection
on $\ol{\mc F}_0$ complete its $(1,0)$-connection to a flat connection.

Both $\on{Op}_{\LG}(X)$ and $\ol{\on{Op}}_{\LG}(X)$ may be viewed as
Lagrangian subspaces of the moduli stack of flat $\LG$-bundles on $X$,
and it turns out that it is their intersection that is relevant to the
eigenfunctions of the global differential operators.

Indeed, we have a large commutative algebra $D_G \otimes \ol{D}_G$ of
global differential operators on the line bundle $K^{1/2} \otimes
\ol{K}^{1/2}$ of half-densities on $\Bun_G$. This algebra is
isomorphic to the algebra of functions on $\on{Op}_{\LG}(X) \times
\ol{\on{Op}}_{\LG}(X)$.

Let $\Bun_G^{\on{st}} \subset \Bun_G$ be the substack of stable
$G$-bundles. Suppose that it is open and dense in $\Bun_G$ (this is
equivalent to the genus of $X$ being greater than 1). We define the
Hilbert space $L^2(\Bun_G)$ as the completion of the space $V$ of
smooth compactly supported sections of $K^{1/2} \otimes \ol{K}^{1/2}$
over $\Bun_G^{\on{st}}$ with the standard Hermitian inner product.

The algebra $D_G \otimes \ol{D}_G$ preserves the space $V$ and is
generated over $\C$ by those operators that are symmetric on
$V$. These are unbounded operators on $L^2(\Bun_G)$, but we expect
that a~real form of the algebra $D_G \otimes \ol{D}_G$ has a canonical
self-adjoint extension (this is explained in \cite{EFK}). If so, then
we get a nice set-up for the problem of finding joint eigenfunctions
and eigenvalues of these operators. It is natural to call these
eigenfunctions the {\em automorphic forms} on $\Bun_G$ (or
$\Bun_G^{\on{st}}$) for a complex algebraic curve. We expect that this
can be generalized to an arbitrary connected reductive complex
group~$G$.

The joint eigenvalues of $D_G \otimes \ol{D}_G$ on $L^2(\Bun_G)$
correspond to points in $\on{Op}_{\LG}(X) \times
\ol{\on{Op}}_{\LG}(X)$, i.e., pairs $(\chi,\rho)$, where $\chi \in
\on{Op}_{\LG}(X)$ and $\rho \in \ol{\on{Op}}_{\LG}(X)$. A joint
eigenfunction corresponding to the pair $(\chi,\rho)$ satisfies the
system of linear PDEs
\begin{gather} \label{PDE}
 H_i \Psi = \chi(H_i) \Psi, \qquad \ol{H}_i \Psi = \rho(\ol{H}_i)
 \Psi,
\end{gather}
where the $H_i$ (resp., the $\ol{H}_i$) are generators of $D_G$ (resp.,
$\ol{D}_G$), and $\chi$ (resp., $\rho$) is viewed as a~homomorphism
$D_G \to \C$ (resp., $\ol{D}_G \to \C$).

\looseness=1 As far as I know, the system \eqref{PDE} was first considered by
Teschner \cite{Teschner}, in the case of $G={\rm SL}_2$ (a similar
idea was also proposed in \cite{F:MSRI}). Teschner did not consider
\eqref{PDE} as a spectral problem in the
sense of self-adjoint operators acting on a Hilbert space. Instead, he
considered the problem of finding the set of real-analytic
single-valued solutions $\Psi$ of the system~\eqref{PDE} in which it
is additionally assumed that $\rho=\ol\chi$ (note that this is not
necessarily so if we do not have a~self-adjointness property). He
outlined in \cite{Teschner} how the solution to this problem can be
related to those ${\rm PGL}_2$-opers (equivalently, projective
connections) $\chi$ on $X$ that have monodromy taking values in the
{\em split real form} ${\rm PGL}_2(\R)$ of ${\rm PGL}_2(\C)$ (up to
conjugation by an element of~${\rm PGL}_2(\C)$).

Projective connections with such monodromy have been described by
Goldman \cite{Goldman}. If the genus of $X$ is greater than 1, then
among them there is a special one, corresponding to the uniformization
of $X$. But there are many other ones as well, and they have been the
subject of interest for many years. It is fascinating that they now
show up in the context of the Langlands correspondence for complex
curves.

In \cite{EFK}, we discuss the spectral problem associated to the
system \eqref{PDE} for a general simply-connected simple Lie group
$G$. Though the $H_i$ and the $\ol{H}_i$ correspond to unbounded
operators on the Hilbert space $L^2(\Bun_G)$, we conjecture that their
linear combinations $(H_i+\ol{H}_i)$ and $(H_i-\ol{H}_i)/{\rm i}$ have
canonical self-adjoint extensions. Furthermore, we conjecture (and
prove in some cases) that the corresponding eigenvalues are the pairs
$(\chi,\rho)$ such that $\rho=\tau(\ol\chi)$, where $\tau$ is the
Chevalley involution, and $\chi$ is an $\LG$-oper on $X$ whose
monodromy representation is isomorphic to its complex conjugate (we
expect that this is equivalent to the property that the monodromy
takes values in the split real form of $\LG$, up to conjugation).

In the next subsection I will illustrate how these opers appear in
the abelian case.

\subsection[The spectra of global differential operators for $G={\rm GL}_1$]{The spectra of global differential operators for $\boldsymbol{G={\rm GL}_1}$} \label{diff GL1}

For simplicity, consider the elliptic curve $X = E_{\rm i} = C/(\Z + \Z{\rm i})$
discussed in Section~\ref{anell}. We identify the neutral component
$\Pic^0(X)$ with~$X$ using a reference point $p_0$, as in Section~\ref{anell}. Then the algebra $D_{{\rm GL}_1}$ (resp.~$\ol{D}_{{\rm GL}_1}$) coincides with the algebra of constant holomorphic
(resp. anti-holomorphic) differential operators on~$X$:
\begin{gather*}
D_{{\rm GL}_1} = \C[\pa_z], \qquad \ol{D}_{{\rm GL}_1} = \C[\pa_{\ol{z}}].
\end{gather*}

The eigenfunctions of these operators are precisely the Fourier
harmonics $f_{m,n}$ given by formula \eqref{fmn}:
\begin{gather*}
f_{m,n} = {\rm e}^{2\pi{\rm i} mx} \cdot {\rm e}^{2\pi{\rm i} ny}, \qquad m,n \in \Z.
\end{gather*}
If we rewrite it in terms of $z$ and $\ol{z}$:
\begin{gather*}
f_{m,n} = {\rm e}^{\pi z (n+{\rm i}m)} \cdot {\rm e}^{-\pi \ol{z}(n-{\rm i}m)},
\end{gather*}
then we find that the eigenvalues of $\pa_z$ and $\pa_{\ol{z}}$ on
$f_{m,n}$ are $\pi(n+{\rm i}m)$ and $-\pi(n-{\rm i}m)$ respectively. Let us recast
these eigenvalues in terms of the corresponding ${\rm GL}_1$-opers.

By definition, a ${\rm GL}_1$-oper is a holomorphic connection on the
trivial line bundle on~$X$ (see \cite[Section~4.5]{F:rev}). The space of
such connections is canonically isomorphic to the space of holomorphic
one-forms on~$X$ which may be written as $-\lambda \,{\rm d}z$, where $\lambda
\in \C$. An element of the space of ${\rm GL}_1$-opers may therefore be
represented as a holomorphic connection on the trivial line bundle,
which together with its $(0,1)$ part $\pa_{\ol{z}}$ yields the flat
connection
\begin{gather} \label{GL1oper}
\nabla = d - \lambda \,{\rm d}z, \qquad \lambda \in \C.
\end{gather}
Under the isomorphism $\on{Spec} D_{{\rm GL}_1} \simeq \on{Op}_{{\rm GL}_1}(X)$,
the oper~\eqref{GL1oper} corresponds to the eigenvalue $\la$ of
$\pa_z$ (this is why we included the sign in \eqref{GL1oper}).

Likewise, an element of the complex conjugate space
$\ol{\on{Op}}_{{\rm GL}_1}(X)$ is an anti-holomorphic connection on the
trivial line bundle, which together with its $(1,0)$ part $\pa_z$
yields the flat connection
\begin{gather} \label{antiGL1oper}
\ol\nabla = d - \mu \,{\rm d}\ol{z}, \qquad \mu \in \C.
\end{gather}
Under the isomorphism $\on{Spec} \ol{D}_{{\rm GL}_1} \simeq
\ol{\on{Op}}_{{\rm GL}_1}(X)$, the oper~\eqref{antiGL1oper} corresponds to
the eigenvalue $\mu$ of~$\pa_{\ol{z}}$.

We have found above that the eigenvalues of $\pa_z$ and $\pa_{\ol{z}}$
on $L^2(\Bun_{{\rm GL}_1})$ are $\pi(n+{\rm i}m)$ and $-\pi(n-{\rm i}m)$, respectively,
where $m,n \in \Z$. The following lemma, which is proved by a direct
computation, links them to ${\rm GL}_1$-opers with monodromy in ${\rm GL}_1(\R)$.

\begin{lem} The connection~\eqref{GL1oper} $($resp.~\eqref{antiGL1oper}$)$ on the trivial line bundle on $E_{\rm i} = \C(\Z + \Z{\rm i})$ has monodromy taking values in the split real form $\R^\times \subset \C^\times$ if and only if $\lambda = \pi(n+{\rm i}m)$ $($resp.~$\mu = -\pi(n-{\rm i}m))$, where $m,n \in \Z$.
\end{lem}

This lemma generalizes in a straightforward fashion to arbitrary curves and arbitrary abelian groups. Namely, the harmonics ${\rm e}^{2\pi{\rm i} \varphi_{\ga}}$, $\ga \in H^1(X,\Lambda^*(T))$, introduced in Section~\ref{torus} are the eigenfunctions of the global differential operators on $\Bun^0_T(X)$. The $\LT$-oper on~$X$ encoding the eigenvalues of the holomorphic differential operators is the holomorphic connection on the trivial $\LT$-bundle on~$X$
\begin{gather*} 
\nabla_\gamma^{\on{hol}} = d - 2\pi{\rm i} \omega_\ga
\end{gather*}
(compare with formula~\eqref{nablaga}). One can show that its monodromy representation takes values in the split real form of~$\LT$, and conversely, these are all the $\LT$-opers on~$X$ that have real monodromy. Thus, the conjectural description of the spectra of global differential operators on~$\Bun_G$ in terms of opers with split real monodromy (see the end of Section~\ref{global}) holds in the abelian case.

Recall that in the abelian case we also have well-defined Hecke operators. It is interesting to note that they commute with the global differential operators and share the same eigenfunctions. Furthermore, the eigenvalues of the Hecke operators may be expressed in terms of the eigenvalues of the global differential operators.

For non-abelian $G$, we consider global differential operators as
substitutes for Hecke operators. We expect that their eigenvalues are
given by the $\LG$-opers satisfying a special condition: namely, their
monodromy representation $\pi_1(X,p_0) \to \LG$ takes values in the
split real form of $\LG$ (up to conjugation). It is natural to view these homomorphisms as
the Langlands parameters of the automorphic forms for curves
over~$\C$. The details will appear in~\cite{EFK}.

\subsection*{Acknowledgments}
The first version of this paper was based on the notes of my talk at
the Sixth Abel Conference held at University of Minnesota in November
2018. I~thank Roberto Alvarenga, Julia Gordon, Ivan Fesenko, David
Kazhdan, and Raven Waller for valuable discussions.

\pdfbookmark[1]{References}{ref}
\LastPageEnding

\end{document}